\documentclass[11pt]{article}
\usepackage [margin=.8in]{geometry}

\usepackage{amsmath}
\usepackage{wasysym}
\usepackage{latexsym}
\usepackage{amsfonts}
\usepackage{mathrsfs}
\usepackage{amssymb}
\usepackage{ifsym}
\usepackage{dsfont}
\usepackage[all]{xy}
\usepackage{authblk}
\usepackage{lipsum}

\newcommand\blfootnote[1]{%
  \begingroup
  \renewcommand\thefootnote{}\footnote{#1}%
  \addtocounter{footnote}{-1}%
  \endgroup
}

\usepackage{amsfonts}
\usepackage{amssymb}

\begin{document}
\title{Set Theoretical Analogues of the Barwise-Schlipf Theorem}
\author{Ali Enayat}

\maketitle

\begin{abstract}
\noindent We prove the following characterizations of nonstandard models of
ZFC (Zermelo-Fraenkel set theory with the axiom of choice) that have an
expansion to a model of $\mathrm{GB}$ (G\"{o}del-Bernays class theory) plus $%
\mathrm{\Delta }_{1}^{1}$\textrm{-}$\mathrm{CA}$\emph{\ }(the scheme of $%
\Delta _{1}^{1}$-Comprehension). In what follows, $\mathcal{M(\alpha )}:=(%
\mathrm{V}(\alpha ),\in )^{\mathcal{M}}$, $\mathbb{L}_{\mathcal{M}}$ is the
set of formulae of the infinitary logic $\mathbb{L}_{\infty ,\omega }$ that
appear in the well-founded part of $\mathcal{M}$, and $\Sigma _{1}^{1}$-$%
\mathrm{AC}$ is the scheme of $\Sigma _{1}^{1}$-Choice.\medskip

\noindent \textbf{Theorem A.~}\textit{The following are equivalent for a nonstandard
model }$\mathcal{M}$ \textit{of }$\mathrm{ZFC}$ \textit{of any cardinality}%
{:}\medskip

\noindent \textbf{(a)} $\mathcal{M(\alpha )}\prec _{\mathbb{L}_{\mathcal{M}}}%
\mathcal{M}$ \textit{for an unbounded collection of }$\alpha \in \mathrm{Ord}%
^{\mathcal{M}}.$\medskip

\noindent \textbf{(b)} $\left( \mathcal{M},{\mathfrak{X}}\right) \models
\mathrm{GB+\Delta }_{1}^{1}$\textrm{-}$\mathrm{CA}$, \textit{where} ${\mathfrak{%
X}}$ \textit{is the family of} $\mathbb{L}_{\mathcal{M}}$-\textit{definable subsets of} $\mathcal{M}$.\medskip

\noindent \textbf{(c)} \textit{There is} $\mathfrak{X}$ \textit{such that} $(%
\mathcal{M},\mathfrak{X})\models \mathrm{GB}+\Delta _{1}^{1}$-$\mathrm{CA}$%
\emph{$\mathsf{.}$}\medskip

\noindent \textbf{Theorem B.~}\textit{The following are equivalent for a \textbf{countable%
} nonstandard model of }\textrm{ZFC}:\medskip

\noindent \textbf{(a)} $\mathcal{M(\alpha )}\prec _{\mathbb{L}_{\mathcal{M}}}%
\mathcal{M}$ \textit{for an unbounded collection of }$\alpha \in \mathrm{Ord}%
^{\mathcal{M}}.$\medskip

\noindent \textbf{(b)}\emph{\ }\textit{There is} $\mathfrak{X}$ \textit{such
that} $($\emph{${\mathcal{M}},{\mathfrak{X}}$}$)\models \mathrm{GB}+\mathrm{%
\Delta }_{1}^{1}$\textrm{-}$\mathrm{CA}+\Sigma _{1}^{1}$-$\mathrm{AC.}$
\end{abstract}

\affil{University of Gothenburg, Gothenburg, Sweden\newline
\texttt{ali.enayat@gu.se}}

\blfootnote {\textit{Acknowledgments}. I am indebted to the indefatigable camaraderie of Jim Schmerl; this paper would not have seen the light of day without our joint recent work \cite%
{enayat-schmerl}. Thanks also to Victoria Gitman, Zachiri McKenzie, and Kameryn Williams for their interest in this work. I am also grateful to the anonymous referee for meticulously examining the paper and weeding out infelicities. The research presented in this paper was partially supported by the National Science Centre, Poland (NCN), grant number 2019/34/A/HS1/00399.}



\begin{center}
\textbf{1.~INTRODUCTION\bigskip }
\end{center}

The point of departure of this paper is Theorem 1.1 below, which
characterizes recursively saturated models of $\mathrm{PA}$ (Peano
Arithmetic) as precisely those nonstandard models of $\mathrm{PA}$ that are
expandable to models of certain subsystems of second order arithmetic. In
what follows, $\mathrm{ACA}_{0}$ is the well-known subsystem of second order
arithmetic whose first order part is $\mathrm{PA}$, $\mathrm{\Delta }_{1}^{1}
$\textrm{-}$\mathrm{CA}$\emph{\ }(respectively $\mathrm{\Sigma }_{1}^{1}$%
\textrm{-}$\mathrm{AC)}$ is the scheme of $\Delta _{1}^{1}$-Comprehension
(respectively $\mathrm{\Sigma }_{1}^{1}$\textrm{-Choice}), and $\mathrm{Def}(%
\mathcal{M})$ is the family of first order definable (parameters allowed) subsets of $\mathcal{%
M}$.\medskip

\noindent \textbf{1.1.~Theorem.~}(Barwise-Schlipf \cite{bs}) \textit{The
following are equivalent for a nonstandard model }\emph{${\mathcal{M}}$} $%
\mathit{of}$\emph{\ $\mathrm{PA}$}\textit{\ of any cardinality}:\emph{%
\medskip }

\noindent \textbf{(a)}\emph{\ }$\mathcal{M}$ \textit{is recursively saturated%
}.\emph{\ \medskip }

\noindent \textbf{(b) }\textit{There is }$\mathfrak{X}$ \textit{such that} $(%
\mathcal{M},\mathfrak{X})\models \mathrm{ACA}_{0}+\Delta _{1}^{1}$-$\mathrm{%
CA}$. \medskip

\noindent \textbf{(c) }$(\mathcal{M},\mathrm{Def}(\mathcal{M}))\models
\mathrm{ACA}_{0}+\Delta _{1}^{1}$-$\mathrm{CA}+\Sigma _{1}^{1}$-$\mathrm{AC}%
. $\medskip

\noindent The argument for $(a)\Rightarrow (c)$ given by Barwise and Schlipf
used the machinery of admissible set theory. Not long after, an elementary
argument was found by Feferman and Stavi (independently), as reported in
Smory\'{n}ski \cite{smo81}. However, the proof presented for $(b)\Rightarrow
(a)$ by Barwise and Schlipf was shown in \cite{enayat-schmerl} to be
impaired by a significant gap, and additionally, a correct proof (using a
technique not available to Barwise and Schlipf) was presented. Now,
prospering from a 45 year hindsight, we can say that the hard part of
Theorem 1.1 is $(b)\Rightarrow (a)$, and the straightforward part is $%
(a)\Rightarrow (c)$ ($(c)\Rightarrow (b)$ is trivial, of course)$.$%
\medskip

An analogue of Theorem 1.1 in the realm of set theory was presented by
Schlipf, as in Theorem 1.2 below, in which $\mathrm{o}(\mathcal{M})$ is the
ordinal height of the well-founded part of $\mathcal{M}$, and $\mathrm{o}(%
\mathrm{HYP}(\mathcal{M}))$ is the ordinal height of $\mathrm{HYP}(\mathcal{M%
}),$ where $\mathrm{HYP}(\mathcal{M})$ is the least admissible structure
over $\mathcal{M}$, as defined in Barwise's definitive text \cite%
{Barwise-book} on admissible set theory. Theorem 1.2 implies the analogue of
Theorem 1.1 for models of $\mathrm{ZF}$ (in which $\mathrm{PA}$ is replaced
by $\mathrm{ZF}$, and $\mathrm{ACA}_{0}$ is replaced by GB), using Schlipf's
characterization of recursive saturation in terms of $\mathrm{o}(\mathrm{HYP}%
(\mathcal{M}))=\omega $.\medskip

\noindent \textbf{1.2.~Theorem.~}(Schlipf \cite{Schlipf-PAMS}) \textit{The
following two conditions are equivalent for a nonstandard model} $\mathcal{M}$ \textit{of} $\mathrm{ZF}$ \textit{of any cardinality}: \medskip

\noindent \textbf{(a) }\textit{There is }$\mathfrak{X}$ \textit{such that} $(%
\mathcal{M},\mathfrak{X})\models \mathrm{GB}+\Delta _{1}^{1}$-$\mathrm{CA}$.
\medskip

\noindent \textbf{(b)}\emph{\ }$\mathrm{o}(\mathcal{M})=\mathrm{o}(\mathrm{%
HYP}(\mathcal{M}))$, \textit{and} $\mathcal{M}$ \textit{satisfies }ZF
\textit{with replacement and separation for formulae involving predicates
for all relations on} $\mathcal{M}$ \textit{that appear in} $\mathrm{HYP}(%
\mathcal{M})$.\medskip

\noindent \textit{Moreover, if} $\mathcal{M}$\ \textit{is a countable
nonstandard model of }\textrm{ZFC}, \textit{then} (a) \textit{and} (b)
\textit{are equivalent to}:\medskip

\noindent \textbf{(c)} \textit{There is} $\mathfrak{X}$ \textit{%
such that} $(\mathcal{M},\mathfrak{X})\models \mathrm{GB}+\Delta _{1}^{1}$%
\textrm{-}$\mathrm{CA}+\Sigma _{1}^{1}$\textrm{-}$\mathrm{AC}$.\medskip

In a different direction, the paper \cite{enayat-condensable} studies the
family of so-called condensable\ models of ZF, a family that includes all
resplendent models of ZF (and in particular, all countable recursively
saturated models of ZF). In the terminology of \cite{enayat-condensable}, a
model $\mathcal{M}\models \mathrm{ZF}$ is \textit{condensable }if $\mathcal{M%
}\cong \mathcal{M(\alpha )}\prec _{\mathbb{L}_{\mathcal{M}}}\mathcal{M}$ for
some \textquotedblleft ordinal\textquotedblright\ $\alpha \in \mathrm{Ord}^{%
\mathcal{M}}$, where $\mathcal{M(\alpha )}:=(\mathrm{V}(\alpha ),\in )^{%
\mathcal{M}}$ and $\mathbb{L}_{\mathcal{M}}$ is the set of formulae of the
infinitary logic $\mathbb{L}_{\infty ,\omega }$ that appear in the
well-founded part of $\mathcal{M}$. The following theorem gives various
characterizations of the notion of condensability (see Section 2 for the
definitions of the technical notions used in the statement of Theorem
1.3).\medskip

\noindent \textbf{1.3.~Theorem.~}\cite{enayat-condensable} \textit{%
The following are equivalent for a \textbf{countable} nonstandard model }$\mathcal{M}$ \textit{of }$%
\mathrm{ZF}$:

\medskip

\noindent \textbf{(a) }$\mathcal{M}$ \textit{is cofinally condensable, i.e.,}
$\mathcal{M}\cong \mathcal{M(\alpha )}\prec _{\mathbb{L}_{\mathcal{M}}}%
\mathcal{M}$ \textit{for an unbounded collection of }$\alpha \in \mathrm{Ord}%
^{\mathcal{M}}$.\textit{\ }\medskip

\noindent \textbf{(b)} $\mathcal{M}$ \textit{is condensable}.\medskip

\noindent \textbf{(c)} \textit{For some nonstandard }$\gamma \in \mathrm{Ord}%
^{\mathcal{M}}$\textit{ and some} $S\subseteq M$, $S$ \textit{is an amenable\
}$\gamma $-\textit{satisfaction class on }$\mathcal{M}.$\medskip

\noindent \textbf{(d)} $\mathcal{M}(\alpha)\prec _{\mathbb{L}_{\mathcal{M}}}\mathcal{M}$ \textit{for an
unbounded collection of }$\alpha \in \mathrm{Ord}^{\mathcal{M}}$.\medskip

\noindent \textbf{(e)} $\mathcal{M}$ \textit{is }$W$-\textit{%
saturated and} $\mathcal{M}\models \mathrm{ZF}(\mathbb{L}_{\mathcal{M}}).$%
\medskip

\noindent \textit{Moreover, without the assumption of countability of} $%
\mathcal{M}$, \textit{the following implications hold}:

\begin{center}
$(a)\Rightarrow (b)\Rightarrow (c)\Rightarrow (d)\Leftrightarrow (e).%
\footnote{%
We suspect that the implication $(b)\Rightarrow (a)$ fails for some
uncountable model of ZF; the implications $(c)\Rightarrow (b)$ and $%
(d)\Rightarrow (c)$ were shown to be irreversible in \cite%
{enayat-condensable}.
\par
{}}$
\end{center}

\noindent The main results of this paper are Theorems A and B below that tie
Theorems 1.2 and 1.3 together. The proofs of these results do not rely on
machinery from admissible set theory, in particular we obtain a new proof,
from first principles, of the equivalence of (a) and (c) of Theorem 1.2 for
a countable nonstandard model $\mathcal{M}$ of \textrm{ZFC}. Note that if $%
\mathcal{M}$ is $\omega $-nonstandard, then condition (a) in Theorems A and
B below is equivalent to recursive saturation of $\mathcal{M}$. \medskip

\noindent \textbf{Theorem A.~}\textit{The following are equivalent for a nonstandard
model }$\mathcal{M}$ \textit{of }$\mathrm{ZFC}$ \textit{of any cardinality}:

\medskip

\noindent \textbf{(a)} $\mathcal{M(\alpha )}\prec _{\mathbb{L}_{\mathcal{M}}}%
\mathcal{M}$ \textit{for an unbounded collection of }$\alpha \in \mathrm{Ord}%
^{\mathcal{M}}.$\medskip

\noindent \textbf{(b)} $\left( \mathcal{M},{\mathfrak{X}}\right) \models
\mathrm{GB+\Delta }_{1}^{1}$\textrm{-}$\mathrm{CA}$, \textit{where} ${\mathfrak{%
X}}$ \textit{is the family of} $\mathbb{L}_{\mathcal{M}}$-\textit{definable subsets of} $\mathcal{M}$.\medskip

\noindent \textbf{(c)} \textit{There is} $\mathfrak{X}$ \textit{such that} $(%
\mathcal{M},\mathfrak{X})\models \mathrm{GB}+\Delta _{1}^{1}$-$\mathrm{CA}$%
\emph{$\mathsf{.}$}\medskip

\noindent \textbf{Theorem B.~}\textit{The following are equivalent for a \textbf{countable%
} nonstandard model of } $\mathrm{ZFC}$:
\medskip

\noindent \textbf{(a)} $\mathcal{M(\alpha )}\prec _{\mathbb{L}_{\mathcal{M}}}%
\mathcal{M}$ \textit{for an unbounded collection of }$\alpha \in \mathrm{Ord}%
^{\mathcal{M}}$.\medskip

\noindent \textbf{(b)}\emph{\ }\textit{There is} $\mathfrak{X}$ \textit{such
that} $($\emph{${\mathcal{M}},{\mathfrak{X}}$}$)\models \mathrm{GB}+\mathrm{%
\Delta }_{1}^{1}$\textrm{-}$\mathrm{CA}+\Sigma _{1}^{1}$-$\mathrm{AC}$%
.\medskip

\noindent We suspect that Theorem A can be strengthened by weakening $%
\mathrm{ZFC}$ to $\mathrm{ZF}$. As explained in Remark 4.4, in Theorem B, $%
\mathrm{ZFC}$ cannot be weakened to $\mathrm{ZF}$, and the assumption of
countability of $\mathcal{M}$ is essential. The proof of Theorem A is
presented in Section 3, and Theorem B is established in Section 4.\textbf{%
\bigskip }

\begin{center}
\textbf{2.~PRELIMINARIES }
\end{center}
\bigskip

In this section we collect the basic definitions, notations, conventions,
and results that will be used in the statements and proofs of our main
results in Sections 3 and 4.\medskip

\noindent \textbf{2.1.~Definition.~}(Models, languages, and theories) Models
will be represented using calligraphic fonts ($\mathcal{M}$, $\mathcal{N}$,
etc.) and their universes will be represented using the corresponding roman
fonts ($M$, $N$, etc.). In the definitions below, $\mathcal{M}$ is a model
of $\mathrm{ZF}$ and $\in ^{\mathcal{M}}$ is the membership relation of $%
\mathcal{M}$. \medskip

\noindent \textbf{(a) }$\mathrm{Ord}^{\mathcal{M}}$ is the class of
\textquotedblleft ordinals\textquotedblright\ of $\mathcal{M}$, i.e., $%
\mathrm{Ord}^{\mathcal{M}}:=\left\{ m\in M:\mathcal{M}\models \mathrm{Ord}%
(m)\right\} ,$ where $\mathrm{Ord}(x)$ expresses \textquotedblleft $x$ is
transitive and is well-ordered by $\in $\textquotedblright . More generally,
given a class $\mathrm{D}$ whose defining formula is $\delta (x)$, $\mathrm{D%
}^{\mathcal{M}}:=\left\{ m\in M:\mathcal{M}\models \delta (m)\right\} .$%
\medskip

\noindent \textbf{(b) }$\mathcal{M}$ is \textit{nonstandard} if $\in ^{%
\mathcal{M}}$ is ill-founded (equivalently: if $(\mathrm{Ord},\in )^{%
\mathcal{M}}$ is ill-founded). $\mathcal{M}$ is $\omega $-\textit{nonstandard%
} if $\left( \omega ,\in \right) ^{\mathcal{M}}$ is ill-founded.\medskip

\noindent \textbf{(c)} For $c\in M$, $\mathrm{Ext}_{\mathcal{M}}(c)$ is the $%
\mathcal{M}$-extension of $c$, i.e., $\mathrm{Ext}_{\mathcal{M}}(c):=\{m\in
M:m\in ^{\mathcal{M}}c\}.$ We say that a subset $X$ of $M$ \textit{is coded
in} $\mathcal{M}$ if there is some $c\in M$ such that $\mathrm{Ext}_{%
\mathcal{M}}(c)=X.$ For $A\subseteq M$, $\mathrm{Cod}_{A}\mathrm{(}\mathcal{%
M)}$ is the collection of sets of the form $A\cap \mathrm{Ext}_{\mathcal{M}%
}(c)$, where $c\in M$. \medskip

\noindent \textbf{(d)} The \textit{well-founded part} of $\mathcal{M}$,
denoted $\mathrm{WF}(\mathcal{M})$, consists of all elements $m$ of $%
\mathcal{M}$ such that there is no infinite sequence $\left\langle
a_{n}:n<\omega \right\rangle $ with $m=a_{0}$ and $a_{n+1}\in ^{\mathcal{M}%
}a_{n}$ for all $n\in \omega .$ Given $m\in M,$ we say that $m$ \textit{is a
nonstandard element of} $\mathcal{M}$\ if $m\notin \mathrm{WF}(\mathcal{M}).$
We denote the submodel of $\mathcal{M}$ whose universe is $\mathrm{WF}(%
\mathcal{M})$ by $\mathcal{WF}(\mathcal{M}).$ It is well-known that if $%
\mathcal{M}$ is a model of $\mathrm{ZF}$, then $\mathcal{WF}(\mathcal{M})$
satisfies $\mathrm{KP}$ (Kripke-Platek set theory) \cite[Chapter II, Theorem
8.4]{Barwise-book}.

\begin{itemize}
\item It is important to bear in mind that we will identify\textit{\ }$%
\mathrm{WF}(\mathcal{M})$ with its transitive collapse.
\end{itemize}

\noindent \textbf{(e)} $\mathrm{o}(\mathcal{M})$ (read as: \textit{the
ordinal of }$\mathcal{M}$) is the supremum of all ordinals in $\mathrm{WF}(%
\mathcal{M}).$\medskip

\noindent \textbf{(f)} Let $\mathcal{L}_{\mathrm{set}}$ be the usual
vocabulary $\{=,\in \}$ of set theory. In this paper we use $\mathbb{L}%
_{\infty ,\omega }$ to denote the infinitary language based on the
vocabulary $\mathcal{L}_{\mathrm{set}}$. Thus $\mathbb{L}_{\infty ,\omega }$
is a set-theoretic language that allows conjunctions and disjunctions of
\textit{sets} (but not proper classes) of formulae, subject to the
restriction that such infinitary formulae have at most finitely many free
variables. Given a set $\Psi $ of formulae, we denote such conjunctions and
disjunctions respectively as $\bigwedge \Psi $ and $\bigvee \Psi $.

\begin{itemize}
\item In the interest of efficiency, we will treat disjunction and universal
quantification as defined notions.
\end{itemize}

\noindent \textbf{(g)} $\mathbb{L}_{\delta ,\omega }$ is the sublanguage of $%
\mathbb{L}_{\infty ,\omega }$ that allows conjunctions and disjunctions of
sets of formulae of cardinality \textit{less than} $\delta .$ Note that $%
\mathbb{L}_{\omega ,\omega }$ is none other than the usual first order
language of set theory, and that in general the language $\mathbb{L}_{\delta
,\omega }$ only uses finite strings of quantifiers (as indicated by the $%
\omega $ in the subscript). \medskip

\noindent \textbf{(h)} We say that $\mathbb{F}$ is a \textit{fragment} of $%
\mathbb{L}_{\infty ,\omega }$ if $\mathbb{F}$ is a set of formulae of $%
\mathbb{L}_{\infty ,\omega }$ that is closed under subformulae, renaming of
free variables, existential quantification, negation, and conjunction.

\begin{itemize}
\item A fragment of $\mathbb{L}_{\infty ,\omega }$ that plays a central role
in this paper is $\mathbb{L}_{\mathcal{M}}:=\mathbb{L}_{\infty ,\omega }\cap
\mathrm{WF}(\mathcal{M)}$. Note that if $M$ is countable, \textbf{\ }$%
\mathbb{L}_{\mathcal{M}}=\mathbb{L}_{\omega _{1},\omega }\cap \mathrm{WF}(%
\mathcal{M)}$.
\end{itemize}

\noindent \textbf{(i) }Given a fragment $\mathbb{F}$ of $\mathbb{L}_{\infty
,\omega }$, and $\mathcal{L}_{\mathrm{set}}$-structures $\mathcal{N}_{1}$
and $\mathcal{N}_{2}$, we write $\mathcal{N}_{1}\prec _{\mathbb{F}}\mathcal{N%
}_{2}$ to indicate that $\mathcal{N}_{1}$ is a submodel of $\mathcal{N}_{2}$%
, and for all $\varphi (x_{1},\cdot \cdot \cdot ,x_{n})\in \mathbb{F}$ and
all tuples $\left( a_{1},\cdot \cdot \cdot ,a_{n}\right) $ from $N_{1}$, we
have:

\begin{center}
$\mathcal{N}_{1}\models \varphi (a_{1},\cdot \cdot \cdot ,a_{n})$ iff $%
\mathcal{N}_{2}\models \varphi (a_{1},\cdot \cdot \cdot ,a_{n})$.
\end{center}

\noindent \textbf{(j)} Given a fragment $\mathbb{F}$ of $\mathbb{L}_{\infty
,\omega }$, $\mathrm{Th}_{\mathbb{F}}(\mathcal{M)}$ is the set of sentences
(closed formulae) of $\mathbb{F}$ that hold in $\mathcal{M}$, and $\mathrm{ZF%
}(\mathbb{F})$ is the natural extension of ZF in which the usual schemes of
separation and collection are extended to the schemes $\mathrm{Sep}(\mathbb{F%
})$ and $\mathrm{Coll}(\mathbb{F})$ so as to allow formulae in $\mathbb{F}$
to be used for \textquotedblleft separating\textquotedblright\ and
\textquotedblleft collecting\textquotedblright\ (respectively).\medskip

\noindent \textbf{(k)} For $\varphi \in \mathbb{L}_{\infty ,\omega }$, the
\textit{depth} of $\varphi ,$ denoted $\mathrm{Depth}(\varphi )$, is the
ordinal defined recursively by the following clauses:\medskip

\noindent (1) $\mathrm{Depth}(\varphi )=0$, if $\varphi $ is an atomic
formula.

\noindent (2) $\mathrm{Depth}(\varphi )=\mathrm{Depth}(\psi )+1,$ if $%
\varphi =\lnot \psi .$

\noindent (3) $\mathrm{Depth}(\varphi )=\mathrm{Depth}(\psi )+1,$ if $%
\varphi =\exists x\ \psi .$

\noindent (4) $\mathrm{Depth}(\varphi )=\sup \{\mathrm{Depth}(\psi )+1:\psi
\in \Psi \}$, if $\varphi =\bigwedge \Psi .$\medskip

\noindent \textbf{(l)} For an ordinal $\alpha $ we use $\mathrm{D(}\alpha ) $
to denote $\{\varphi \in \mathbb{L}_{\infty ,\omega }:\mathrm{Depth}(\varphi
)<\alpha \}$. Within KP, one can code each formula $\varphi \in $ $\mathbb{L}%
_{\infty ,\omega }$ with a set $\ulcorner \varphi \urcorner $ as in Chapter
3 of \cite{Barwise-book}, but in the interest of better readability we will
often identify a formula with its code. This coding allows us to construe
statements such as $\varphi \in \mathbb{L}_{\infty ,\omega }$ and $\mathrm{%
Depth}(\varphi )=\alpha $ as statements in the first order language of set
theory. It is easy to see that for a sufficiently large $k\in \omega $, $%
\mathrm{D(}\alpha )\subseteq \mathrm{V(}\omega +k\alpha )$ for each ordinal $%
\alpha $. This makes it clear that $\mathbb{L}_{\mathcal{M}%
}=\bigcup\limits_{\alpha \in \mathrm{o}(\mathcal{M})}\mathrm{D}^{\mathcal{M}}%
\mathrm{(}\alpha ).$

\noindent \textbf{(m) }Suppose $\mathcal{M}$ is nonstandard and $W:=\mathrm{%
WF}(\mathcal{M)}$. $\mathcal{M}$ is $W$-\textit{saturated} if for every $%
k\in \omega $ and every type $p(x,y_{1},\cdot \cdot \cdot ,y_{k})$, and for
every $k$-tuple $\overline{a}$ of parameters from $\mathcal{M}$, $p(x,%
\overline{a})$ is realized in $\mathcal{M}$ provided the following three
conditions are satisfied:$\medskip $

\noindent $(m1)$ $p(x,\overline{y})\subseteq \mathbb{L}_{\mathcal{M}}$.$%
\medskip $

\noindent $(m2)$ $p(x,\overline{y})\in \mathrm{Cod}_{W}(\mathcal{M})$.$%
\medskip $

\noindent $(m3)$ $\forall w\in W\ \mathcal{M}\models \exists x\left(
\bigwedge\limits_{\varphi \in p(x,\overline{y})\cap w}\varphi (x,\overline{a}%
)\right) .$

\noindent \textbf{(n) }Every model of $\mathrm{GB}$ can be put in the form $%
\left( \mathcal{N},\mathfrak{X}\right) ,$ where $\mathcal{N}\models \mathrm{%
ZF}$ and $\mathfrak{X}\subseteq \mathcal{P}(N).$\medskip

\noindent \textbf{2.2.~Definition.~}Suppose $\mathcal{M}$ is a model of $%
\mathrm{ZF}$, and $S\subseteq M$. \medskip

\noindent \textbf{(a)} $S$ is \textit{separative} (over $\mathcal{M}$) if $(%
\mathcal{M},S)$ satisfies the separation scheme $\mathrm{Sep(S)}$ in the
extended language that includes a fresh predicate \textrm{S} (interpreted by
$S$). \medskip

\noindent \textbf{(b)} $S$ is \textit{collective} (over $\mathcal{M}$) if $(%
\mathcal{M},S)$ satisfies the collection scheme $\mathrm{Coll(S)}$ in the
extended language that includes a fresh predicate \textrm{S} (interpreted by
$S$).\medskip

\noindent \textbf{(c)} $S$ is \textit{amenable} (over $\mathcal{M}$) if $S$
is both separative and collective$.$ In other words, $S$ is amenable if $%
\left( \mathcal{M},S\right) $ satisfies the replacement scheme $\mathrm{Repl}%
(S)$ in the extended language that includes a fresh predicate \textrm{S}
(interpreted by $S$). Note that if $\left( \mathcal{M},\mathfrak{X}\right) $
is a model of GB, then each element of $\mathfrak{X}$ is amenable over $%
\mathcal{M}$.\medskip

\noindent \textbf{(d)} For $\alpha \in \mathrm{Ord}^{\mathcal{M}}$, $S$ is
an $\alpha $\textit{-satisfaction class} (over $\mathcal{M}$) if $S$
correctly decides the truth of atomic sentences, and $S$ satisfies Tarski's
compositional clauses of a truth predicate for \textrm{D}$^{\mathcal{M}}$%
\textrm{(}$\alpha )$-sentences (see below for the precise definition). $S$
is an $\infty $\textit{-satisfaction class} over $\mathcal{M}$ if $S$ is an
$\alpha $\textit{-}satisfaction class over $\mathcal{M}$ for every $\alpha
\in \mathrm{Ord}^{\mathcal{M}}$. \medskip

\noindent We elaborate the meaning of (d) above. Reasoning within ZF, for
each object $a$ in the universe of sets, let $c_{a\text{ }}$ be a constant
symbol denoting $a$ (where the map $a\mapsto c_{a}$ is $\Delta _{1}),$ and
let $\mathrm{Sent}^{+}(\alpha ,x)$ be the set-theoretic formula (with an
ordinal parameter $\alpha $ and the free variable $x)$ that defines the
proper class of sentences of the form $\varphi \left( c_{a_{1}},\cdot \cdot
\cdot ,c_{a_{n}}\right) $, where $\varphi (x_{1},\cdot \cdot \cdot
,x_{n})\in \mathrm{D}$\textrm{(}$\alpha )$ (the superscript $+$ on $\mathrm{%
Sent}^{+}(\alpha ,x)$ indicates that $x$ is a sentence in the language
augmented with the indicated proper class of constant symbols). Then $S$ is
an $\alpha $-satisfaction class over $\mathcal{M}$ if $\left( \mathcal{M}%
,S\right) \models \mathrm{Sat}(\alpha ,S)$, where $\mathrm{Sat}(\alpha ,S)$
is the conjunction of the universal generalizations of axioms $(I)$
through $(IV)$ below:\medskip

\noindent $(I)\ \  \left( \mathrm{S}\left( \ulcorner
c_{a}=c_{b}\urcorner \right) \leftrightarrow a=b\right) \wedge \left(
\mathrm{S}\left( \ulcorner c_{a}\in c_{b}\urcorner \right) \leftrightarrow
a\in b\right)  .$\medskip

\noindent $(II)\ \ \left( \mathrm{Sent}^{+}(\alpha ,\varphi )\wedge \left(
\varphi =\lnot \psi \right) \right) \rightarrow \left( \mathrm{S}(\varphi
)\leftrightarrow \lnot \mathrm{S}\mathsf{(}\psi \mathsf{)}\right) \mathsf{.}$%
\medskip

\noindent $(III)$ $\ \left( \mathrm{Sent}^{+}(\alpha ,\varphi )\wedge \left(
\varphi =\bigwedge \Psi \right) \right) \rightarrow \left( \mathrm{S}%
(\varphi )\leftrightarrow \forall \psi \in \Psi \ \mathrm{S}\mathsf{(}\psi
\mathsf{)}\right) \mathsf{.}$\medskip

\noindent $(IV)$ $\ \left( \mathrm{Sent}^{+}(\alpha ,\varphi )\wedge \left(
\varphi =\exists x\ \psi (x)\right) \right) \rightarrow \left( \mathrm{S}%
(\varphi )\leftrightarrow \exists x\ \mathrm{S}\mathsf{(\psi (}c_{x}\mathsf{%
))}\right) .$\medskip

\noindent \textbf{(e)} For $\alpha <\mathrm{o}(\mathcal{M})$, $S$ is \textit{%
the} $\alpha $-satisfaction class over $\mathcal{M},$ if $S$ is the usual
Tarskian satisfaction class for formulae in $\mathbb{L}_{\mathcal{M}}$ of
depth less than $\alpha ,$ i.e., the unique $\alpha $-satisfaction class $S$
over $\mathcal{M}$ such that $S$ satisfies:\medskip

\noindent $(V)$ $\ \forall x\left( \mathrm{S}(x)\rightarrow \mathrm{Sent}%
^{+}(\alpha ,x)\right) .$\medskip

\noindent Finally, the $\mathrm{o}(\mathcal{M})$-\textit{satisfaction class
over} $\mathcal{M}$ is the usual Tarskian satisfaction class for formulae in
$\mathbb{L}_{\mathcal{M}}$ of depth less than $\mathrm{o}(\mathcal{M}),$
i.e., the union of all Tarskian $\alpha $-satisfaction classes over $%
\mathcal{M}$ as $\alpha $ ranges in $\mathrm{o}(\mathcal{M}).$\medskip

\begin{itemize}
\item In the interest of a lighter notation, if $S$ is an $\alpha $\textit{%
-satisfaction class} over $\mathcal{M}$ (for a possibly nonstandard $\alpha
\in \mathrm{Ord}^{\mathcal{M}})$, $\varphi (x_{1},\cdot \cdot \cdot ,x_{n})$
is an $n$-ary formula of \textrm{D}$^{\mathcal{M}}$\textrm{(}$\alpha )$, and
$a_{1},\cdot \cdot \cdot ,a_{n}$ are in $M,$ we will often write $\varphi
\left( a_{1},\cdot \cdot \cdot ,a_{n}\right) \in S$ instead of $\varphi
\left( c_{a_{1}},\cdot \cdot \cdot ,c_{a_{n}}\right) \in S.$
\end{itemize}

The following proposition is immediately derivable from the relevant
definitions.\medskip

\noindent \textbf{2.3.~Proposition.~}\textit{\ If} $S$ \textit{is an} $\alpha
$\textit{-satisfaction class} \textit{over} $\mathcal{M}$ \textit{for some
nonstandard ordinal} $\alpha $ \textit{of} $\mathcal{M}$, \textit{then for
all }$n$\textit{-ary formula }$\varphi (x_{1},\cdot \cdot \cdot ,x_{n})$
\textit{of} $\mathbb{L}_{\mathcal{M}}$ \textit{and all }$n$-\textit{tuples} $%
(a_{1},\cdot \cdot \cdot ,a_{n})$ \textit{from} $M$, \textit{we have}:

\begin{center}
$\mathcal{M}\models \varphi (a_{1},\cdot \cdot \cdot ,a_{n})$ \textit{iff} $\varphi
(a_{1},\cdot \cdot \cdot ,a_{n})\in S.$
\end{center}

\noindent \textit{In particular, for all sentences} $\varphi $ \textit{of} $%
\mathbb{L}_{\mathcal{M}},$ $\varphi \in S$ \textit{iff} $\varphi \in \mathrm{%
Th}_{\mathbb{L}_{\mathcal{M}}}(\mathcal{M)}.$\medskip

\noindent \textbf{2.4.~Remark.~}Reasoning within ZF, given any limit ordinal
$\gamma ,$ $\left( \mathrm{V}(\gamma),\in \right) $ carries a separative $\gamma $%
-satisfaction class $S$ since we can take $S$ to be the Tarskian
satisfaction class on $\left( \mathrm{V}(\gamma),\in \right) $ for formulae of
depth less than $\gamma .$ More specifically, the Tarski recursive
construction/definition of truth works equally well in this more general
context of infinitary languages since $\left( \mathrm{V}(\gamma),\in \right) $
forms a set. Observe that $\left( \mathrm{V}(\gamma),\in ,S\right) \models \mathrm{%
Sep(S)}$ comes \textquotedblleft for free\textquotedblright\ since for any $%
X\subseteq \mathrm{V}(\gamma)$ the expansion $\left( \mathrm{V}(\gamma),\in ,X\right) $
satisfies the scheme of separation in the extended language. \medskip

\noindent \textbf{2.5.~}\textbf{Proposition.~}\cite[Proposition 2.5]%
{enayat-condensable}\textbf{\ }\textit{If}$\mathcal{M}\models \mathrm{KP}$,
\textit{then for each nonzero }$\alpha \in \mathrm{o}(\mathcal{M})$ \textit{there is
a formula} $\mathrm{Sat}_{\alpha }(x)\in \mathbb{L}_{\mathcal{M}}$ \textit{%
such that }$\mathrm{Sat}_{\alpha }^{\mathcal{M}}(x)$ \textit{is the} $\alpha
$-\textit{satisfaction class over} $\mathcal{M}.$ \medskip

The following general version of the elementary chains theorem of model
theory can verified by a routine adaptation of
the usual proof of the $\mathbb{L}_{\omega ,\omega }$-version (e.g., as in
\cite[Theorem 3.1.9]{Chang-Keisler}).  \medskip

\noindent \textbf{2.6.~}\textbf{Proposition.~}(Elementary Chains) \textit{%
Suppose }$\mathbb{L}\subseteq \mathbb{L}_{\infty ,\omega }$ \textit{where }$%
\mathbb{L}$ \textit{is closed under subformulae;} $(I,\vartriangleleft )$
\textit{is a linear order; }$\left\langle \mathcal{M}_{i}:i\in
I\right\rangle $ \textit{is an} $\mathbb{L}$-\textit{elementary chain }(%
\textit{i.e.,} $\mathcal{M}_{i}\prec _{\mathbb{L}}\mathcal{M}_{j}$ \textit{%
whenever} $i\vartriangleleft j$); \textit{and} $\mathcal{M}%
=\bigcup\limits_{i\in I}\mathcal{M}_{i}.$\textit{\ Then }$\mathcal{M}%
_{i}\prec _{\mathbb{L}}\mathcal{M}$ \textit{for each} $i\in I.$\medskip

\medskip

The following generalization of the Montague-Vaught reflection theorem of set theory appears as Proposition 2.8 of \cite{enayat-condensable}, in a slightly weaker form, where the class of $\varphi$-reflecting ordinals (where $\varphi$ ranges over $\mathrm{D}^{\mathcal{M}}(\alpha )$) is asserted to be unbounded, as opposed to closed and unbounded. The stronger version below can be readily obtained by putting Proposition 2.6 above together with Proposition 2.8 of \cite{enayat-condensable}. The closed unboundedness of the class of $\varphi$-reflective ordinals is needed in the proof of Lemma 4.3 of this paper, where it is important to arrange arbitrarily large $\varphi$-reflective ordinals of \textit{countable cofinality}. \medskip

\noindent \textbf{2.7.~Proposition.~}(Reflection)\textbf{\ }\textit{Suppose }%
$\mathcal{M}\models \mathrm{ZF}(\mathbb{L}_{\mathcal{M}}),$ \textit{and for
each} $\varphi \in \mathbb{L}_{\mathcal{M}}$ \textit{where} $\varphi $
\textit{is} $n$\textit{-ary, let }$\mathrm{Ref}_{\varphi }(\gamma )$ \textit{%
be the }$\mathbb{L}_{\mathcal{M}}$-\textit{formula}:

\begin{center}
$\forall x_{1}\in \mathrm{V}(\gamma )\cdot \cdot \cdot \forall x_{n}\in
\mathrm{V}(\gamma )\ \left[ \varphi \left( x_{1},\cdot \cdot \cdot
,x_{n}\right) \longleftrightarrow \varphi ^{\mathrm{V}(\gamma )}\left(
x_{1},\cdot \cdot \cdot ,x_{n}\right) \right]. $
\end{center}

\noindent \textit{Then for any} $\alpha \in \mathrm{o}(\mathcal{M})$ \textit{%
there is a closed unbounded collection of ordinals} $\gamma \in \mathrm{Ord}^{\mathcal{M}}$ \textit{%
such that }$\mathcal{M}(\gamma )$ \textit{reflects all formulae in} $\mathrm{D}^{\mathcal{M}}$($\alpha ),$ \textit{i.e}., $\mathcal{M}\models \mathrm{Ref}%
_{\varphi }(\gamma )$ \textit{for all }$\mathbb{L}_{\mathcal{M}}$\textit{%
-formulae} $\varphi $ \textit{of depth less than} $\alpha$. \medskip

The notion \textquotedblleft $f$ is $\lambda $-onto $X$\textquotedblright\
introduced in Definition 2.8 below, and the corresponding existence result
(Proposition 2.9), are adaptations of Definition 2.1 and Lemma 2.2 of
Kaufmann and Schmerl's work \cite{KaufmannSchmerlAPAL} on models of
arithmetic to the setting of set theory. Lemma 2.2 of \cite%
{KaufmannSchmerlAPAL} played a key role in the proof presented in \cite%
{enayat-schmerl} of the direction $(b)\Rightarrow (a)$ of Theorem 1.1.
Proposition 2.9 plays an analogous role in the proof of the direction $%
(b)\Rightarrow (a)$ of Theorem A.\medskip

\noindent \textbf{2.8.~Definition.~}A set $I$ is an \textit{ordinal interval}
if $I=\{\gamma :\alpha <\gamma <\beta \}$ for some ordinals $\alpha $ and $%
\beta .$ Suppose $f:I\rightarrow X$, where $I$\ is an ordinal interval, $X$
is some set, and $\lambda $ is an ordinal. The notion $f$ \textit{is} $%
\lambda $-\textit{onto} $X$ is defined by recursion on $\lambda $ as follows:

\begin{itemize}
\item $f$ is $0$-onto\ $X$ means: $f$ is onto $X$.

\item For $\lambda =\gamma +1$, $f$ is $\lambda $-onto $X$ means: for each $%
Y\subseteq X$ there is an ordinal interval $J\subseteq I$ such that $%
f\upharpoonright J$ is $\gamma $-onto $Y$.

\item For a limit ordinal $\lambda ,f$ is $\lambda $-onto\ $X$ means: $%
\forall \gamma <\lambda $ $f$ is $\gamma $-onto $X$.\medskip
\end{itemize}

\noindent \textbf{2.9.~Proposition.~} (ZFC) \textit{Given any set }$X$ \textit{%
and any ordinal }$\lambda $ \textit{there is some ordinal interval} $I$%
\textit{\ and a function} $f:I\rightarrow X$ \textit{such that} $f$ \textit{%
is} $\lambda $-\textit{onto }$X$\textit{.}\medskip

\noindent \textbf{Proof.~}We use induction on $\lambda .$ The case $\lambda
=0$ is clear since we are working in $\mathrm{ZFC}$. If $\lambda =\gamma +1$, then
for each subset $Y$ of $X$ there is an ordinal interval $I_{Y}$ and some function $h_{Y}:I_{Y}\rightarrow Y$ such that $h_Y$ is $\gamma$-onto $Y$. Use AC to enumerate $%
\mathcal{P}(X)$ as $\left\{ Y_{\alpha }:\alpha <\kappa \right\} $, where $%
\kappa =\left\vert \mathcal{P}(X)\right\vert .$ Let $I_{\alpha
}:=I_{Y_{\alpha }}$ and $h_{\alpha}:=h_{Y_{\alpha}}$ for each $\alpha <\kappa $ and choose an ordinal
interval $I$ that is order isomorphic to the well-ordering $%
\sum\limits_{\alpha <\kappa }I_{\alpha }$. More explicitly, let $Z=\left\{
\{\alpha \}\times I_{\alpha }:\alpha <\kappa \right\} $ and let $%
\vartriangleleft $ be the lexicographic order on $Z.$ Then since $%
\vartriangleleft $ is a well-ordering, there is an ordinal interval $I$ and
an isomorphism $F$ between $(Z,\vartriangleleft )$ and $\left( I,\in \right)
.$ Note that $F(\{\alpha \}\times I_{\alpha })\cap F(\{\beta \}\times
I_{\beta })=\varnothing $ when $\alpha $ and $\beta $ are distinct elements
of $\kappa$. Since for each $\alpha< \kappa$ the function $h_{\alpha }:I_{\alpha }\rightarrow
Y_{\alpha} $ has the property of being $\lambda $-onto $Y_{\alpha}$, the isomorphism $F$
allows us to construct functions $f_{\alpha }:F(\{\alpha \}\times I_{\alpha
})\rightarrow Y_{\alpha}$ such that each $f_{\alpha }$ is $\lambda $-onto $Y_{\alpha}.$ This
will ensure that $\bigcup\limits_{\alpha <\kappa }f_{\alpha }$
is a function from $I$ to $X$ that is $\gamma +1$-onto $X$.\medskip

For limit $\lambda $ we use a strategy similar to the successor case. By
inductive assumption, for each $\gamma <\lambda $ there is some ordinal
interval $I_{\gamma }$ and a function $f_{\gamma }:$
$I_{\gamma }\rightarrow X$ such that $f_{\gamma }$ is $\gamma $-onto $X$. We can
therefore find an ordinal interval $I$ and an isomorphism $F$ between the
well-ordering $\sum\limits_{\gamma <\lambda }I_{\gamma }$ and $I$. For each $%
\gamma <\lambda $ we can then construct functions $f_{\gamma }:F(\{\gamma
\}\times I_{\gamma })\rightarrow X$ such that $f_{\gamma }$ is $\gamma $%
-onto $X$. It is evident that $\bigcup\limits_{\alpha <\lambda }f_{\alpha }$
is a function from $I$ to $X$ that is $\lambda$-onto $X$. \hfill $\square $

\medskip

\noindent \textbf{2.10.~Remark.~}Recall that $\Sigma^{1}_{k}$-$\mathrm{AC}$
(AC for the axiom of choice) is the scheme consisting of formulae of the
form

\begin{center}
$\forall x\ \exists X\ \psi (x,X)\rightarrow \exists Y\ \forall x\ \psi
(x,\left( Y\right) _{x}),$
\end{center}

\noindent where $\psi (x,X)$ is a $\Sigma^{1}_{k}$-formula (parameters allowed), and $\Sigma^{1}_{k}$-$\mathrm{Coll}$ (Coll for Collection) is the
scheme consisting of formulae of the form

\begin{center}
$\forall x\ \exists X\ \psi (x,X)\rightarrow \exists Y\ \forall x\ \exists
y\ \psi (x,\left( Y\right) _{y}),$
\end{center}

\noindent where $\psi (x,X)$ is a $\Sigma^{1}_{k}$-formula (again, with parameters allowed). In the above

\begin{center}
$\left( Y\right) _{x}:=\left\{ y:\left\langle x,y\right\rangle \in Y\right\}
$,
\end{center}

\noindent where $\left\langle x,y\right\rangle $ is a canonical pairing
function. \medskip

\noindent \textbf{(a)} It is well-known that $\Sigma _{k}^{1}$-$\mathrm{AC}$
implies $\Delta _{k}^{1}$-$\mathrm{CA}$ for all $k<\omega $; an easy proof
in the arithmetical setting can be found in \cite[Lemma~VII.6.6(1)]{sim};
the same proof readily works in the set-theoretic context. \medskip

\noindent \textbf{(b)} Let $\mathrm{GBC}$ be the result of augmenting $%
\mathrm{GB}$ with the global axiom of choice. It is well-known that in the
presence of $\mathrm{GBC,}$ (1) $\Sigma _{1}^{1}$-$\mathrm{AC}$ is
equivalent to $\Sigma _{1}^{1}$-$\mathrm{Coll}$, and (2) global choice is
provable in $\mathrm{GB}+\Sigma _{1}^{1}$-$\mathrm{AC.}$ For more detail,
see, e.g., \cite[Section 3.1]{Fujimoto APAL}.\medskip

\noindent Let $\overline{\mathbb{L}}_{\infty ,\omega }$ be the extension of $%
\mathbb{L}_{\infty ,\omega }\mathcal{\ }$based on the extended vocabulary $%
\mathcal{L}_{\mathrm{set}}=\{=,\in ,f\}$, where $f$ is a function symbol
(for a global choice function), and let $\overline{\mathbb{L}}_{\mathcal{N}%
}:=\overline{\mathbb{L}}_{\infty ,\omega }\cap \mathrm{WF}(\mathcal{N)}$,
where $\mathcal{N}\models \mathrm{ZF}$\ ($\mathcal{N}$\ need not be
nonstandard, so $\mathrm{WF}(\mathcal{N)}$ might be the whole of $\mathcal{N}
$). The following result is the infinitary generalization of the well-known
theorem that global choice can be generically added to models of ZFC of
countable cofinality \cite{Felgner} (and its proof is similar to the proof of the finitary case). A proof of part (b) of Proposition 2.11 can
be found in \cite[Theorem 11]{Schlipf-PAMS}. \medskip

\noindent \textbf{2.11.~Proposition.~}(Forcing Global Choice)\textbf{\ }\textit{%
Let }$\mathcal{N}\models \mathrm{ZFC}(\mathbb{L}_{\mathcal{N}}),$ \textit{and%
} $\mathbb{P}$\textit{\ be the class notion of forcing consisting of set
choice functions in} $\mathcal{N}$, \textit{ordered by set inclusion}.
\medskip

\noindent \textbf{(a)} \textit{If} $\mathrm{Ord}^{\mathcal{N}}$ \textit{has
countable cofinality, then there is an} $\mathbb{L}_{\mathcal{N}}$\textit{%
-generic filter} $G\subseteq \mathbb{P}$\textit{, in the sense that }$G$%
\textit{\ is a filter that intersects every dense subset of} $\mathbb{P}$
\textit{that is definable in }$\mathcal{N}$ \textit{by a formula in} $%
\mathbb{L}_{\mathcal{N}}$ (\textit{parameters allowed}).\medskip

\noindent \textbf{(b)} \textit{If }$G$ \textit{is} \textit{an }$\mathbb{L}_{%
\mathcal{N}}$-\textit{generic filter over }$\mathbb{P}$\textit{, and }$%
f=\cup G,$ \textit{then }$f$ \textit{is a global choice function over} $%
\mathcal{N}$,\textit{\ and} $(\mathcal{N},f)\models \mathrm{ZF}(\overline{%
\mathbb{L}}_{\mathcal{N}}).$

\bigskip
\pagebreak
\begin{center}
\textbf{3.~PROOF OF THEOREM A}\bigskip
\end{center}

In this section we establish the first main result of our paper. In part (b) of the
following theorem, $\mathrm{Def}_{\mathbb{L}_{\mathcal{M}}}$ is the family of $\mathbb{L}_{\mathcal{M}}$-definable subsets of $M$ (parameters allowed).

\medskip

\noindent \textbf{3.1.~Theorem.~}\textit{The following are equivalent for a
nonstandard model }$\mathcal{M}$ \textit{of} $\mathrm{ZF}$ \textit{of any
cardinality}:\medskip

\noindent \textbf{(a)} $\mathcal{M(\alpha )}\prec _{\mathbb{L}_{\mathcal{M}}}%
\mathcal{M}$ \textit{for an unbounded collection of }$\alpha \in \mathrm{Ord}%
^{\mathcal{M}}.$\medskip

\noindent \textbf{(b)} $\left( \mathcal{M},{\mathfrak{X}}\right) \models
\mathrm{GB+\Delta }_{1}^{1}$\textrm{-}$\mathrm{CA}$, for $\mathfrak{X} = \mathrm{Def}_{\mathbb{L}_{\mathcal{M}}}$.

\medskip

\noindent \textbf{(c)} \textit{There is} $\mathfrak{X}$ \textit{such that} $(%
\mathcal{M},\mathfrak{X})\models \mathrm{GB}+\Delta _{1}^{1}$-$\mathrm{CA}$%
\emph{$\mathsf{.}$}\medskip

\noindent \textbf{Proof.~}Since $(b)\Rightarrow (c)$ is trivial, it suffices
to establish $(a)\Rightarrow (b)$ and $(c)\Rightarrow (a).$\medskip

\noindent $\mathbf{(a)\Rightarrow (b)}.$ Assume (a). Then by $%
(d)\Leftrightarrow (e)$ of Theorem 1.3, we have:\medskip

\noindent (1) \emph{${\mathcal{M}}$} satisfies $\mathrm{ZF}$\emph{$\mathsf{(}
$}$\mathbb{L}_{\mathcal{M}})$, and\medskip

\noindent (2) \emph{${\mathcal{M}}$} is $W$-saturated.\medskip\

\noindent (1) makes it clear that $\mathrm{GB}$ holds in $\left( \mathcal{M},%
{\mathfrak{X}}\right) $. We will use (2) to show that \emph{$\Delta _{1}^{1}$%
-$\mathrm{CA}$ }holds in $\left( \mathcal{M},{\mathfrak{X}}\right) .$ To
this end, let $U\subseteq M$ such that $U$ is defined in $\left( \mathcal{M},%
{\mathfrak{X}}\right) $ by a $\Sigma _{1}^{1}$-formula $\exists X\ \psi
^{+}(X,x,A)$, and $M\backslash U$ is defined in $\left( \mathcal{M},{%
\mathfrak{X}}\right) $ by a $\Sigma _{1}^{1}$-formula $\exists X\ \psi
^{-}(X,A,x)$, where $A\in {\mathfrak{X}}$ is a class parameter definable by
the $\mathbb{L}_{\mathcal{M}}$-formula $\alpha (m,v)$. Here $m\in M$ is a
set parameter; note that we may assume without loss of generality that the
only parameter in $\psi ^{+}$ and in $\psi ^{-}$ is a class parameter $A$.
Consider the infinitary formulae $\theta ^{+}(x)$ and $\theta ^{-}(x)$
defined as follows:\medskip

\begin{center}
$\theta ^{+}(x):=\bigvee\limits_{\varphi (y,v)\in \mathbb{L}_{\mathcal{M}%
}}\exists y$\ $\psi ^{+}(X/\varphi (y,v),A/\alpha (m,v),x)$; and \medskip

$\theta ^{-}(x):=\bigvee\limits_{\varphi (y,v)\in \mathbb{L}_{\mathcal{M}%
}}\exists y$\ $\psi ^{-}(X/\varphi (y,v),A/\alpha (m,v),x).$
\end{center}

\noindent In the above $\psi ^{+}(X/\varphi (y,v),A/\alpha (m,v),x)$
(respectively $\psi ^{-}(X/\varphi (y,v),A/\alpha (m,v),x)$) is the result
of replacing all occurrences of subformulae of the form $w\in X$ (where $w$
is a variable) in $\psi ^{+}$ (respectively in $\psi ^{-}$) by $\varphi
(y,w),$ and replacing all occurrences of subformulae of the form $w\in A$ in
$\psi ^{+}$ (respectively in $\psi ^{-}$) by $\alpha (m,w)$, and re-naming
variables to avoid unintended clashes. Since each $X\in {\mathfrak{X}}$ can
be written in the form $\left\{ v\in M:\mathcal{M}\models \varphi
(m_{1},v)\right\} $ (where $m_{1}\in M$ is a parameter), $U$ is definable in
$\mathcal{M}$\ by $\theta ^{+}(x)$ and $M\backslash U$ is definable in $%
\mathcal{M}$\ by $\theta ^{-}(x).$ Therefore we have:\medskip

\noindent (3) $\mathcal{M}\models \forall x\left( \theta ^{+}(x)\vee \theta
^{-}(x)\right) .$\medskip

\noindent Next, we aim to verify (4) below. In what follows $\mathrm{D}^{\mathcal{M}}(\alpha)$ is as in part (k) of Definition 2.1.\medskip

\noindent (4) There is some $\alpha \in \mathrm{o}(\mathcal{M})$ such that $%
\mathcal{M}\models \forall x\left( \theta _{\alpha }^{+}(x)\vee \theta
_{\alpha }^{-}(x)\right) $, where

\begin{center}
$\theta _{\alpha }^{+}(x):=\bigvee\limits_{\varphi (y,v)\in \mathrm{D}^{\mathcal{M}}(\alpha)}\exists y\ \psi ^{+}(X/\varphi(y,v),A/\alpha (m,v),x)$; and\medskip

$\theta _{\alpha }^{-}(x):=\bigvee\limits_{\varphi (y,v)\in \mathrm{D}^{\mathcal{M}}(\alpha)}\exists y\ \psi ^{-}(X/\varphi(y,v),A/\alpha (m,v),x).$
\end{center}

\noindent Notice that (4) implies that $U$ is definable in $\mathcal{M}$ by $%
\theta _{\alpha }^{+}(x)$, so the verification of $\Delta _{1}^{1}$-$\mathrm{%
CA}$ will be complete once we establish (4) since $\theta _{\alpha
}^{+}(x)\in \mathbb{L}_{\mathcal{M}}$ and $\left\{ a\in M:\mathcal{M}\models
\theta _{\alpha }^{+}(a)\right\} \in \mathfrak{X.}$ To establish (4) we
argue by contradiction. Suppose \medskip

\noindent (5) $\mathcal{M}\models \exists x\ \lnot \left( \theta _{\alpha
}^{+}(x)\vee \theta _{\alpha }^{-}(x)\right) $ for each $\alpha \in \mathrm{o%
}(\mathcal{M})$.\medskip

\noindent Consider the $\mathbb{L}_{\mathcal{M}}$-type $p(x),$where

\begin{center}
$p(x):=\left\{ \lnot \left( \theta _{\alpha }^{+}(x)\vee \theta _{\alpha
}^{-}(x)\right) :\alpha \in \mathrm{o}(\mathcal{M})\right\} .$
\end{center}

\noindent It is easy to see that $p(x)\in \mathrm{Cod}_{W}(\mathcal{M}).$ By
(5), for each $\alpha \in \mathrm{o}(\mathcal{M})$, $p(x)\cap \mathrm{D}^{\mathcal{M}}
(\alpha )$ is realized in $\mathcal{M}$, so by $W$-saturation of $\mathcal{M}
$, $p(x)$\ is realized in $\mathcal{M}$, i.e., $\mathcal{M}\models \exists
x\ \lnot \left( \theta ^{+}(x)\vee \theta ^{-}(x)\right) ,$ which
contradicts (3) and finishes the proof of (4)\medskip

\noindent $\mathbf{(c)\Rightarrow (a)}.$ This is the hard direction of
Theorem 3.1 and will require a good deal of preliminary lemmata. It will be
proved as Lemma 3.6. In part (a) of Lemma 3.2, \textrm{Sat}$_{\alpha }$ is
as in Proposition 2.5.\medskip

\noindent \textbf{3.2.~Lemma.~}\textit{If}\emph{\ }$\left( \mathcal{M},{%
\mathfrak{X}}\right) \models \mathrm{GB}+\Delta _{1}^{1}$-$\mathrm{CA}$,%
\textit{\ then the following hold}:\medskip

\noindent \textbf{(a)} \textrm{Sat}$_{\alpha }^{\mathcal{M}}\in {\mathfrak{X}%
}$ \textit{for each nonzero} $\alpha \in \mathrm{o}(\mathcal{M}).$\medskip

\noindent \textbf{(b)} $\mathcal{M}\models \mathrm{ZF}(\mathbb{L}_{\mathcal{M%
}}).$\medskip

\noindent \textbf{Proof.~}To see that (a) holds we will use induction on $%
\alpha $ to verify that \textrm{Sat}$_{\alpha }^{\mathcal{M}}$ is $\Delta
_{1}^{1}$-definable in $($\emph{${\mathcal{M}},{\mathfrak{X}}$}$)$\emph{\ }%
for each $\alpha \in \mathrm{o}(\mathcal{M})$. Recall that $\mathrm{Sat}%
(\alpha ,S)$ is the first order formula that expresses \textquotedblleft $S$
is an $\alpha $-satisfaction class\textquotedblright\ (as in Definition
2.2). Suppose \textrm{Sat}$_{\alpha }^{\mathcal{M}}\in {\mathfrak{X}}$ for
some $\alpha \in \mathrm{o}(\mathcal{M}){\mathfrak{.}}$ Then for each $m\in
M $ we have:

\begin{center}
$m\in \mathrm{Sa}$\textrm{t}$_{\alpha +1}^{\mathcal{M}}$ iff \smallskip

$(\mathcal{M},{\mathfrak{X}})\models \exists S\left[ \mathrm{Sat}(\alpha
,S)\wedge \left( \mathrm{Depth}(m)<\alpha \right) \wedge \left( \mathrm{Neg}%
(m)\vee \mathrm{Exist}(m)\vee \mathrm{Conj}(m)\right) \right] $,
where\smallskip

$\mathrm{Neg}(x):=\exists y\left[ \left( x=\ulcorner \lnot y\urcorner
\right) \wedge \lnot S(y)\right] ;$\smallskip

$\mathrm{Exist}(x):=\exists y\ \exists v\left[ \left( x=\ulcorner \exists v\
y(v)\urcorner \right) \wedge \exists v\ S(y(c_{v}))\right] ;$ and\smallskip

$\mathrm{Conj}(x):=\exists y\left[ \left( x=\ulcorner \bigwedge y\urcorner
\right) \wedge \forall z\in y\ S(z)\right] .$\smallskip
\end{center}

\noindent Similarly, for each $m\in M$ we have:

\begin{center}
$m\in \mathrm{Sa}$\textrm{t}$_{\alpha +1}^{\mathcal{M}}$ iff \smallskip

$(\mathcal{M},{\mathfrak{X}})\models \forall S\left[  \left(\mathrm{Sat}%
(\alpha ,S)\wedge  \mathrm{Depth}(m)=\alpha \right)
\rightarrow \left( \mathrm{Neg}(m)\vee \mathrm{Exist}(m)\vee \mathrm{Conj}%
(m)\right) \right] .$
\end{center}

\noindent Thus $\mathrm{Sa}$\textrm{t}$_{\alpha +1}^{\mathcal{M}}$ has both
a $\Sigma _{1}^{1}$ and a $\Pi _{1}^{1}$ definition in $(\mathcal{M},{%
\mathfrak{X}}).$ The limit case is more straightforward since for limit $%
\alpha $ the following hold for each $m\in M:$\medskip

\begin{center}
$m\in \mathrm{Sa}$\textrm{t}$_{\alpha }^{\mathcal{M}}$ iff $(\mathcal{M},{%
\mathfrak{X}})\models \exists \beta <\alpha \ \left[ \left( \mathrm{Depth}%
(m)=\beta \right) \wedge \exists S\left( \mathrm{Sat}(\beta ,S)\wedge
S(m)\right) \right] ,$ and\medskip

$m\in \mathrm{Sa}$\textrm{t}$_{\alpha }^{\mathcal{M}}$ iff $(\mathcal{M},{%
\mathfrak{X}})\models \exists \beta <\alpha \ \left[ \left( \mathrm{Depth}%
(m)=\beta \right) \wedge \forall S\left( \mathrm{Sat}(\beta ,S)\rightarrow
S(m)\right) \right] .$
\end{center}

\noindent This concludes the proof of (a). Note that (b) is an immediate
consequence of (a) since the veracity of any $\mathbb{L}_{\mathcal{M}}$%
-instance of replacement in $\mathcal{M}$ follows from the amenability of
\textrm{Sat}$_{\alpha }^{\mathcal{M}}$ over $\mathcal{M}$ for a sufficiently
large\textit{\ }$\alpha \in \mathrm{o}(\mathcal{M}).$\hfill $\square $
(Lemma 3.2)\medskip

The notion of paradefinability introduced in Definition 3.3 below is the
set-theoretical analogue\ of the notion of recursive $\sigma $-definability
in \cite{enayat-schmerl}.\medskip

\noindent \textbf{3.3.~Definition.~}Suppose ${\mathcal{M}}\models \mathrm{KP}
$ and $A\subseteq M$. Then, $A$ is \textit{paradefinable} in $\mathcal{M}$,
if there is a sequence $\langle \varphi _{\alpha }(x,\overline{y}):\alpha <%
\mathrm{o}(\mathcal{M})\rangle \in \mathrm{Cod}_{W}(\mathcal{M})$ of $%
\mathbb{L}_{\mathcal{M}}$-formulae (where $\overline{y}$ is a finite tuple
whose length is independent of $\alpha $) such that for some fixed tuple of
parameters $\overline{m}$ in $M$ of the same length as $\overline{y}$, each $\varphi _{\alpha }(x,\overline{m})$
defines a subset $A_{\alpha }\subseteq M$ (in $\mathcal{M)}$, with $%
A=\bigcup\limits_{\alpha <\mathrm{o}(\mathcal{M})}A_{\alpha }$. Under these
conditions, we say that $A$ is \textit{paradefinable} \textit{by} $\langle
\varphi _{\alpha }(x,\overline{m}):\alpha <\mathrm{o}(\mathcal{M})\rangle $.

\medskip

\noindent \textbf{3.4.~Lemma.~}If ${\mathcal{M}}\models \mathrm{KP}$,
\textit{then the following are} \textit{paradefinable in} $\mathcal{M}$%
:\medskip

\noindent \textbf{(a)} $\mathrm{o}(\mathcal{M}).$\medskip

\noindent \textbf{(b) }$\mathrm{WF}(\mathcal{M})$.\medskip

\noindent \textbf{(c)} \textit{The} $\mathrm{o}(\mathcal{M})$-\textit{%
satisfaction class on} $\mathcal{M}$.\medskip

\noindent \textbf{Proof.~}

\noindent \textbf{(a)} $\mathrm{o}(\mathcal{M})$ is paradefinable in $%
\mathcal{M}$ by $\langle E_{\alpha }(x):\alpha <\mathrm{o}(\mathcal{M}%
)\rangle $, where $E_{\alpha }(x)$ (which defines $\left\{ \alpha \right\} )$
is constructed by recursion via:

\begin{center}
$E_{0}(x):=\forall y(y\notin x),$ and for $\alpha >0$, $E_{\alpha
}(x):=\forall y\left( y\in x\leftrightarrow \bigvee\limits_{\lambda <\alpha
}E_{\lambda }(y)\right) .$
\end{center}

\noindent \textbf{(b)\ }$\mathrm{WF}(\mathcal{M})$ is paradefinable in $%
\mathcal{M}$ by $\langle \exists y\left( E_{\alpha }(y)\wedge x\in \mathrm{V}%
(y)\right) :\alpha <\mathrm{o}(\mathcal{M})\rangle .$\medskip

\noindent \textbf{(c)} The $\mathrm{o}(\mathcal{M})$-satisfaction class on $%
\mathcal{M}$ is paradefinable in $\mathcal{M}$ by $\langle \mathrm{Sat}%
_{\alpha }(x):\alpha <\mathrm{o}(\mathcal{M})\rangle $.\hfill $\square $
(Lemma 3.4)\medskip

Part (b)\ of the next lemma is the set-theoretical analogue of \cite[Lemma 1
(b)]{enayat-schmerl}.\medskip

\noindent \textbf{3.5.~Lemma.~}\textit{Suppose}\emph{\ }$\left( \mathcal{M},{%
\mathfrak{X}}\right) \models \mathrm{GB}+\Delta _{1}^{1}$-$\mathrm{CA}$,
\textit{for some nonstandard} $\mathcal{M}$ \textit{that} \textit{is not} $W$%
-\textit{saturated. Then}:\textit{\medskip }

\noindent \textbf{(a)} \textit{For all} $\delta \in \mathrm{Ord}^{\mathcal{M}%
}$ $\left( \delta \in \mathrm{o}(\mathcal{M})\Longleftrightarrow \exists
S\in {\mathfrak{X\ }}\left( \mathcal{M},S\right) \models \mathrm{Sat}(\delta
,S)\right) .$\textit{\medskip }

\noindent \textbf{(b)} \textit{If} $A\subseteq M$ \textit{is} \textit{%
paradefinable in }$\mathcal{M}$\textit{, then} $A$ \textit{is} $\Sigma
_{1}^{1}$-\textit{definable in} $(\mathcal{M},\mathfrak{X})$.\medskip

\noindent \textbf{Proof.~}To establish (a),\textbf{\ }first note
that the assumption of the failure of $W$-saturation in $\mathcal{M}$ by Theorem
1.3 implies that there is no $S\in {\mathfrak{X}}$ such that $S$\ is a $%
\gamma $-satisfaction class over $\mathcal{M}$ for any nonstandard $\gamma
\in \mathrm{Ord}^{\mathcal{M}}$. Combined with part (a) of Lemma 3.2, this
makes it clear that (a) holds. To verify (b), let $A$ be paradefinable by $%
\langle \varphi _{\alpha }(x,\overline{m}):\alpha <\mathrm{o}(\mathcal{M}%
)\rangle $. By replacing $\varphi _{\alpha }(x,\overline{y})$ with $\bigvee\limits_{\beta
\leq \alpha }\varphi _{\beta }(x,\overline{y})$, we can assume that $\mathrm{Depth}%
(\varphi _{\alpha }(x,\overline{y}))<\mathrm{Depth}(\varphi _{\beta }(x,%
\overline{y}))$ for all $\alpha <\beta <\mathrm{o}(\mathcal{M})$. Let $%
\delta $ be a nonstandard element of $\mathrm{Ord}^{\mathcal{M}}$ such that $%
\langle \varphi _{\alpha }(x,\overline{y}):\alpha <\mathrm{\delta }\rangle $
is in $\mathcal{M}$ and extends $\langle \varphi _{\alpha }(x,\overline{y}%
):\alpha <\mathrm{o}(\mathcal{M})\rangle $ and $\mathrm{Depth}(\varphi
_{\alpha }(x,\overline{y}))<\mathrm{o}(\mathcal{M})$ for all $\alpha <%
\mathrm{o}(\mathcal{M})$. Then $A$ is $\Sigma _{1}^{1}$-definable in $({%
\mathcal{M}},{\mathfrak{X}})$ by the formula $\exists X\ \theta (x,X)$,
where
\begin{equation*}
\theta (x,X)=\exists \gamma \lbrack \mathrm{Sat}(\gamma ,X)\wedge \exists
\alpha <\delta \text{ }\mathrm{Depth}(\varphi _{\alpha }(x,\overline{y}%
))<\gamma \wedge \varphi _{\alpha }(c_{x},\overline{m})\in X].
\end{equation*}%
By part (a) of the lemma, this makes it evident that $A$ is $\Sigma _{1}^{1}$%
-definable in $({\mathcal{M}},{\mathfrak{X}})$. \hfill $\square $ (Lemma
3.5)\medskip

\begin{itemize}
\item The proof of Theorem 3.1 will be complete once we verify Lemma 3.6
below, which takes care of the direction $\mathbf{(c)\Rightarrow (a)}$ of
Theorem 3.1. The proof of Lemma 3.6 is rather complicated and we therefore beg the reader's indulgence.
\end{itemize}

\noindent \textbf{3.6.~Lemma.~} \textit{If}\emph{\ }$\left( \mathcal{M},{%
\mathfrak{X}}\right) \models \mathrm{GB}+\Delta _{1}^{1}$-$\mathrm{CA}$
\textit{and} \emph{${\mathcal{M}}$}\textit{\ is nonstandard, then}\emph{\ }$%
\mathcal{M(\alpha )}\prec _{\mathbb{L}_{\mathcal{M}}}\mathcal{M}$ \textit{%
for an unbounded collection of\ }$\alpha \in \mathrm{Ord}^{\mathcal{M}}.$%
\medskip

\noindent \textbf{Proof.~} Suppose not, then by Theorem 1.3, \emph{${\mathcal{%
M}}$} is not $W$-saturated, so by part (a)\ of Lemma 3.5 we conclude:\medskip

\noindent (1) For all $\delta \in \mathrm{Ord}^{\mathcal{M}}$ $\left( \delta
\in \mathrm{o}(\mathcal{M})\Longleftrightarrow \exists S\in {\mathfrak{X\ }}%
\left( \mathcal{M},S\right) \models \mathrm{Sat}(\delta ,S)\right) .$\medskip

\noindent By our supposition there is some $\beta \in \mathrm{Ord}^{\mathcal{M}}$ such that:\medskip

\noindent (2) There is no $\alpha \in \mathrm{Ord}^{\mathcal{M}}$ above $%
\beta $ with $\mathcal{M}(\alpha )\prec _{\mathbb{L}_{\mathcal{M}}}\mathcal{M%
}.$ \medskip

\noindent Since $\mathcal{M}\models \mathrm{ZF}(\mathbb{L}_{\mathcal{M}})$,
by Theorem 2.7 (Reflection), in the real world there is a sequence $%
\left\langle \gamma _{\alpha }:\alpha \in \mathrm{o}(\mathcal{M}%
)\right\rangle $ of ordinals of $\mathcal{M}$ such that for each $\alpha \in
\mathrm{o}(\mathcal{M})$ the following holds:\smallskip

\begin{center}
$\mathcal{M}\models $ \textquotedblleft $\gamma _{\alpha }$ is the first
ordinal $\gamma >\beta $ such that $\mathrm{V}(\gamma )\prec _{\mathrm{D}%
(\alpha )}\mathrm{V}$\textquotedblright $,$
\end{center}

\noindent where $\mathrm{V}(x)\prec _{\mathrm{D}(\alpha )}\mathrm{V}$
abbreviates the $\mathbb{L}_{\mathcal{M}}$-formula $\bigwedge\limits_{\varphi \in \mathrm{D}(\alpha )}%
\mathrm{Ref}_{\varphi }(x)$ ($\mathrm{Ref}_{\varphi }$ is as in Proposition 2.7).

\noindent Let $\Gamma =\{\gamma _{\alpha }:\alpha \in \mathrm{o}(\mathcal{M}%
)\}.$ Clearly $\gamma _{\alpha }\leq \gamma _{\xi }$ for $\alpha <\xi <%
\mathrm{o}(\mathcal{M}).$

\begin{itemize}
\item We now distinguish between the following two cases, and will show that
each leads to a contradiction, thus proving Lemma 3.6 (recall that the proof of Lemma 3.6 starts with ``Suppose not''). Our proof was
inspired by the proof of \cite[Theorem 4]{enayat-schmerl}, which the reader
is highly advised to review before reading the proof below, especially since
the argument for Case B below is a more complex version of the argument for
the \textquotedblleft tall case\textquotedblright\ in the proof of \cite[%
Theorem 4]{enayat-schmerl}. One of the reasons for this increase in
complexity has to do with the fact that in nonstandard models of arithmetic
(and in $\omega $-nonstandard models of set theory) it is easy to find an
ill-founded subset $A$ of the nonstandard ordinals of $\mathcal{M}$ that is
paradefinable since we can choose $A$ to be $\{c-n:n\in \omega \}$, where $c$
is any nonstandard finite ordinal. The existence of such an ill-founded $A$
plays a key role in the proof of \cite[Lemma 2.4]{KaufmannSchmerlAPAL} since
in conjunction with the arithmetical analogue of Proposition 2.9, it allows one
to deduce that if some recursive type is omitted, then a recursive type
consisting of formulae describing an ordinal interval is omitted. However,
in a nonstandard model $\mathcal{M}$ of set theory that is $\omega $%
-standard, the existence of an ill-founded paradefinable subset of the
nonstandard ordinals of $\mathcal{M}$ takes far more effort to establish (with the help of additional assumptions, as indicated in the proof of Case B below).
\end{itemize}

\noindent \textbf{Case A:~} $\Gamma $\textit{\ }is cofinal in $\mathrm{Ord}^{%
\mathcal{M}}.$ We wish to show that $\Gamma $ is $\Delta _{1}^{1}$-definable
in $\left( \mathcal{M},\mathfrak{X}\right) .$ By part (b) of Lemma 3.5, it
is sufficient to show that both $\Gamma $ and its complement are
paradefinable in $\mathcal{M}$. The definition of $\Gamma $ makes it clear
that $\Gamma $ is paradefinable\textit{\ }in $\mathcal{M}$\textit{\ }by $%
\left\langle \varphi _{\alpha }^{+}(x,\beta ):\alpha \in \mathrm{o}(\mathcal{%
M)}\right\rangle ,$ where $\varphi _{\alpha }^{+}(x,\beta )$ is the
following formula:

\begin{center}
$\left( x\in \mathrm{Ord}\wedge \beta \in x\right) \wedge \left( \mathrm{V}%
(x)\prec _{\mathrm{D}(\alpha )}\mathrm{V}\right) \wedge \forall y\in x{%
\mathfrak{\ }}\lnot \left( \mathrm{V}(y)\prec _{\mathrm{D}(\alpha )}\mathrm{V%
}\right) .$
\end{center}

\noindent To see that the complement of $\Gamma $ is also paradefinable in $%
\mathcal{M},$ observe that $\gamma _{\alpha }\leq \gamma _{\xi }$ whenever
$\alpha \leq \xi <\mathrm{o}(\mathcal{M})$, and by Proposition 2.6
(Elementary Chains) $\Gamma $ is a closed subset of $\mathrm{Ord}^{\mathcal{M%
}},$ i.e., for limit $\xi \in \mathrm{o}(\mathcal{M})$, $\gamma _{\xi
}=\sup \{\gamma _{\alpha }:\alpha <\xi \}.$ Thus for each $\nu \in \mathrm{%
Ord}^{\mathcal{M}}\backslash \Gamma ,$ there is some $\alpha \in \mathrm{o}(%
\mathcal{M})$ such that $\gamma _{\alpha }<\nu <\gamma _{\alpha +1}.$ So the
complement of $\Gamma $ is paradefinable in $\mathcal{M}$ by $\left\langle
\theta (x)\vee \varphi _{\alpha }^{-}(x,\beta ):\alpha \in \mathrm{o}(%
\mathcal{M)}\right\rangle ,$ where:

\begin{center}
$\theta (x):=x\notin \mathrm{Ord,}$\medskip

$\varphi _{0}^{-}(x,\beta ):=\left[ \exists y\left( \varphi _{0}^{+}(y,\beta
)\wedge x\in y\right) \right] $, and\medskip

for $\alpha >0,\ \varphi _{\alpha }^{-}(x,\beta ):=\left[ \exists y\ \exists
z\left( \varphi _{\alpha }^{+}(y,\beta )\wedge \varphi _{\alpha
+1}^{+}(z,\beta )\wedge (y\in x\in z)\right) \right] .$\medskip
\end{center}

\noindent Therefore $\Gamma \in {\mathfrak{X}}$, which implies that $\Gamma $
is amenable over $\mathcal{M}$, so coupled with the fact that $\Gamma $ is
cofinal in $\mathrm{Ord}^{\mathcal{M}}$ we conclude that there is some $f\in
\mathfrak{X}$ such that $f$ is an isomorphism between $\Gamma $ and $\mathrm{%
Ord}^{\mathcal{M}}$ (both ordered by $\in ^{\mathcal{M}})$. This contradicts
the fact that $\Gamma $ is well-founded and $\mathrm{Ord}^{\mathcal{M}}$ is
ill-founded, and thus shows that Case A is impossible.\medskip

\noindent \textbf{Case B:~}$\Gamma $ is bounded in $\mathrm{Ord}^{\mathcal{M}%
}.$ In this case, by (2) the supremum of $\Gamma $ does not exist in $%
\mathrm{Ord}^{\mathcal{M}}.$ Fix an upper bound $\delta \in \mathrm{Ord}^{%
\mathcal{M}}$ for $\Gamma $, and $\mathrm{f}$or each $\alpha \in \mathrm{o}(%
\mathcal{M})$ let

\begin{center}
$\psi _{\alpha }(x,\beta ,\delta ):=\left( \beta \in x \in \delta
\right) \wedge \left( \mathrm{V}(x)\prec _{\mathrm{D}(\alpha )}\mathrm{V}%
\right) .$
\end{center}

\noindent Note that if $\alpha $ and $\xi $ are in $\mathrm{o}(\mathcal{M}%
) $ with $\alpha \leq \xi $, then $\mathcal{M}\models \forall x\left( \psi
_{\beta }(x,\xi ,\delta )\rightarrow \psi _{\alpha }(x,\xi ,\delta
)\right) .$ Consider the $\mathbb{L}_{\mathcal{M}}$-type:

\begin{center}
$p(x,\beta ,\delta ):=\left\{ \psi _{\alpha }(x,\beta ,\delta ):\alpha \in
\mathrm{o}(\mathcal{M})\right\} .$
\end{center}

\noindent Clearly:\medskip

\noindent (3) $p(x,\overline{y})\in \mathrm{Cod}_{W}(\mathcal{M}),$ $%
\mathcal{M}\models \exists x\ \psi _{\alpha }(x,\beta ,\delta )$ for each $%
\alpha \in \mathrm{o}(\mathcal{M})$.\medskip

\noindent Moreover, (2) implies:\medskip

\noindent (4) $\mathcal{M}\models \forall x\bigvee\limits_{\alpha \in
\mathrm{o}(\mathcal{M})}\lnot \psi _{\alpha }(x,\beta ,\delta ).$

\noindent In the real world define $\left\langle \delta _{\alpha }:\alpha
\in \mathrm{o}(\mathcal{M})\right\rangle $ with:

\begin{center}
$\delta _{0}=\delta $, and $\delta _{\alpha }=\max \left\{ \xi \in \delta :%
\mathrm{V}(\xi )\prec _{\mathrm{D}(\alpha )}\mathrm{V}\right\} $.
\end{center}

\noindent It is easy to see that $\delta _{\alpha }$s are
well-defined for each $\alpha \in \mathrm{o}(\mathcal{M})$. More
specifically, let $X_{\alpha }=\left\{ \xi \in \delta :\mathrm{V}(\xi )\prec
_{\mathrm{D}(\alpha )}\mathrm{V}\right\} .$ Then by the choice of $\delta ,$
$X_{\alpha }$ is nonempty for each $\alpha \in \mathrm{o}(\mathcal{M})$; and
by part (b) of Lemma 3.2, $X_{\alpha }$ is coded in $\mathcal{M}$, so $\sup
\left( X_{\alpha }\right) $ is well-defined, and by Proposition 2.6
(Elementary Chains) $\sup \left( X_{\alpha }\right) \in X_{\alpha }$, so $%
\max (X_{\alpha })$ is well-defined. It should also be clear that:\medskip

\noindent (5) $\left\{ \delta _{\alpha }:\alpha \in \mathrm{o}(\mathcal{M}%
)\right\} $ is paradefinable in $\mathcal{M}$,\medskip

\noindent (6) $\delta _{\alpha }\geq \delta _{\beta }$ if $\alpha \leq \beta
\in \mathrm{o}(\mathcal{M}),$ and $\medskip $

\noindent (7) $\delta _{\alpha }>\gamma _{\nu }>\beta $ if $\nu \in \mathrm{o%
}(\mathcal{M})$ and $\alpha \in \mathrm{o}(\mathcal{M}).$ \medskip

\noindent Next, we observe:\medskip

\noindent (8) For each $\alpha \in \mathrm{o}(\mathcal{M})$ $\exists \beta
\in \mathrm{o}(\mathcal{M})$ such that $\beta >\alpha $ and $\delta _{\alpha
}>\delta _{\beta }$. \medskip

\noindent To see that (8) is true, note that $\delta _{\alpha }>\beta $ by
(7), so if (8) is false, then $\mathrm{V}(\delta _{\alpha })\prec _{\mathrm{D%
}(\beta )}\mathrm{V}$ for all $\beta \in \mathrm{o}(\mathcal{M})$, which
contradicts (2). Thus (8) implies that $\left\{ \delta _{\alpha }:\alpha \in
\mathrm{o}(\mathcal{M})\right\} $ is ill-founded when viewed as a subset of
\textrm{Ord}$^{\mathcal{M}}.$ Moreover, for any $\alpha \in \mathrm{o}(%
\mathcal{M}),$ $\left\{ \delta _{\beta }:\beta <\alpha \right\} $ is finite.
To verify this, first note that there is a fixed natural number $k$ such
that the depth of the $\mathbb{L}_{\mathcal{M}}$-formula that defines $%
\delta _{\beta }$ for any $\beta \in \mathrm{o}(\mathcal{M})$ is at most $%
\beta +k.$ Therefore in light of (5) and (6) and the fact that $\mathrm{Sat}%
_{\alpha }^{\mathcal{M}}$ (which is present in $\mathfrak{X}$ by part (a) of
Lemma 3.2. and is therefore amenable over $\mathcal{M}$) can evaluate the
defining formulae of $\left\{ \delta _{\beta }:\beta <\alpha \right\} ,$ the
well-foundedness of $\mathrm{Ord}^{\mathcal{M}}$ as viewed in $\mathcal{M}$
implies that $\left\{ \delta _{\beta }:\beta <\alpha \right\} $ is finite
from the point of view of $\mathcal{M}$. Therefore $\left\{ \delta _{\beta
}:\beta <\alpha \right\} $ is finite in the real world (this is clear if $%
\mathcal{M}$ is $\omega $-standard; if $\mathcal{M}$ is $\omega $%
-nonstandard, then it is trivial since $\alpha $ would have to be a finite
ordinal since $\alpha \in \mathrm{o}(\mathcal{M})).$ Putting all this
together, we conclude:\medskip

\noindent (9) The order type of $\left\{ \delta _{\alpha }:\alpha \in
\mathrm{o}(\mathcal{M})\right\} $ under $\in ^{\mathcal{M}}$ is $\omega
^{\ast }$ (i.e., the reversal of $\omega $).\medskip

We will now use the sequence $\left\langle \delta _{\alpha }:\alpha \in
\mathrm{o}(\mathcal{M})\right\rangle $ to describe a type $\widehat{p}%
(x,\beta ,\delta )$ such that $\widehat{p}(x,\overline{y})\in \mathrm{Cod}%
_{W}(\mathcal{M})$, with $\widehat{p}(x,\overline{y})=\left\{ \widehat{\psi }%
_{\alpha }(x,\overline{y}):\alpha \in \mathrm{o}(\mathcal{M})\right\} $,
where each $\widehat{\psi }_{\alpha }(x,\beta ,\delta )$ describes an
ordinal interval, i.e.,

\begin{center}
$\widehat{\psi }_{\alpha }(x,\beta ,\delta ):=s_{\alpha }(\beta ,\delta
)<x<t_{\alpha }(\beta ,\delta ),$
\end{center}

\noindent for an appropriate choice of $\mathbb{L}_{\mathcal{M}}$-definable
terms $\left\langle s_{\alpha }(\beta ,\delta ):\alpha \in \mathrm{o}(%
\mathcal{M})\right\rangle $ and $\left\langle t_{\alpha }(\beta ,\delta
):\alpha \in \mathrm{o}(\mathcal{M})\right\rangle $, where each $s_{\alpha
}(\beta ,\delta )$ and $t_{\alpha }(\beta ,\delta )$ is in $\mathrm{Ord}^{%
\mathcal{M}},$ and $s_{\alpha }$ and $t_{\alpha }$ are defined below.\medskip

Let $\mathrm{I}(x,y):=\{z:x\in z\in y\}.$ In $\mathcal{M}$, define $X$ as
the ordinal interval $\mathrm{I}(\beta ,\delta )$, and apply Proposition 2.9 to
get hold of a function $f$ and some ordinal interval $I$ such that $%
f:I\rightarrow X$ and $f$ is $\delta $-onto $X.$ Let $I_{0}$ be the $%
\vartriangleleft $-first ordinal interval $I$ that supports such a function,
where $\vartriangleleft $ is a canonical well-ordering of all ordinal
subintervals. Then define $s_{0}$ and $t_{0}$ so that $I_{0}=\mathrm{I}%
(s_{0}(\beta ,\delta ),t_{0}(\beta ,\delta ))$. For $\alpha >0$ we define $%
s_{\alpha }$ and $t_{\alpha }$ by recursion on $\alpha :$\medskip

\begin{itemize}
\item If $\alpha $ is a successor ordinal $\lambda +1$, then $s_{\alpha
}(\beta ,\delta )$ and $t_{\alpha }(\beta ,\delta )$ are respectively the
left and right end points of the $\vartriangleleft $-first ordinal
subinterval $I$ of the ordinal interval $\mathrm{I}(s_{\lambda }(\beta
,\delta ),t_{\lambda }(\beta ,\delta ))$ such that $f\upharpoonright I$ is $%
\delta _{\alpha }$-onto $\{x\in \mathrm{I}(\beta ,\delta ):\psi _{\alpha
}(x,\beta ,\delta )\}\mathrm{.}$

\item If $\alpha $ is a limit ordinal, then $s_{\alpha }(\beta ,\delta )$
and $t_{\alpha }(\beta ,\delta )$ are respectively the left and right end
points of the $\vartriangleleft $-first (ordinal) subinterval $I$ of \textrm{%
I}$(s_{\lambda _{0}}(\beta ,\delta ),t_{\lambda _{0}}(\beta ,\delta ))$ such
that $f\upharpoonright I$ is $\delta _{\alpha }$-onto $\{x:\psi _{\alpha
}(x,\beta ,\delta )\}$, where $\lambda _{0}$ is the first ordinal below $%
\alpha $ for which the tail $\left\langle \delta _{\lambda }:\lambda
_{0}\leq \lambda <\alpha \right\rangle $ is a constant sequence.
\end{itemize}

\noindent Next we will show: \medskip

\noindent (10) $\mathcal{M}\models \exists x\ \widehat{\psi }_{\alpha
}(x,\beta ,\delta )$ for each $\alpha \in \mathrm{o}(\mathcal{M}).$ \medskip

\noindent Naturally, we use induction on $\alpha $ to verify (10). Proposition
2.9 and part (b) of Lemma 3.2 make it clear that the induction smoothly goes
through for the base case and the successor case. The limit case requires
the additional fact that if $\alpha $ is a limit ordinal, then by (9) there is some $%
\lambda _{0}<\alpha $ such that the tail $\left\langle \delta _{\lambda
}:\lambda _{0}\leq \lambda <\alpha \right\rangle $ is a constant sequence.
\medskip

\noindent Finally we will establish: \medskip

\noindent (11) $\widehat{p}(x,\beta ,\delta )$ is not realized in $\mathcal{M%
}$.\medskip

\noindent To verify (11) recall that within $\mathcal{M}$, $f$ maps each
interval $\mathrm{I}(s_{\alpha }(\beta ,\delta ),t_{\alpha }(\beta ,\delta
)) $ into $\left\{ x:\psi _{\alpha }(x,\beta ,\delta )\right\} $. Therefore
if some element $m$ of $\mathcal{M}$ realizes $\widehat{p}(x,\beta ,\delta )$%
, then $f(m)$ realizes $p(x,\beta ,\delta ),$ which contradicts (4). We are
now finally ready to wrap up the proof. Let

\begin{center}
$I:=\left\{ x\in \mathrm{Ord}^{\mathcal{M}}:\exists \alpha \in \mathrm{o}(%
\mathcal{M})\left( x<s_{\alpha }\left( \beta ,\delta \right) \right)
\right\} .$
\end{center}

\noindent It is evident that $I$ is paradefinable in $\mathcal{M}$. The
complement of $I$ can written as:

\begin{center}
$M\backslash I=\left\{ x:x\notin \mathrm{Ord}^{\mathcal{M}}\vee \exists
\alpha \in \mathrm{o}(\mathcal{M})\left( x>t_{\alpha }\left( \beta ,\delta
\right) \right) \right\} ,$
\end{center}

\noindent which makes it clear that $M\backslash I$ is also paradefinable in
$\mathcal{M}$. Therefore by part (b) of Lemma 3.5 both $I$ and its
complement are $\Sigma _{1}^{1}$-definable in $(\mathcal{M},\mathfrak{X})$
and thus $I\in \mathfrak{X}$, which implies that the supremum of $I$ exists
in $\mathcal{M}$ (since each element of $\mathfrak{X}$ is separative over $%
\mathcal{M)}$. This contradicts (11) and concludes the demonstration that
Case B is impossible.\hfill $\square $ (Lemma 3.6 and Theorem 3.1)\bigskip

\begin{center}
\textbf{4.~PROOF OF THEOREM B}\bigskip
\end{center}

In this section we establish the second main result of this paper.\bigskip

\noindent \textbf{4.1.~Theorem.~} \textit{The following are equivalent for a
countable nonstandard model }$\mathcal{M}$ \textit{of} $\mathrm{ZFC}$%
:\medskip

\noindent \textbf{(a)} $\mathcal{M(\alpha )}\prec _{\mathbb{L}_{\mathcal{M}}}%
\mathcal{M}$ \textit{for an unbounded collection of }$\alpha \in \mathrm{Ord}%
^{\mathcal{M}}.$\medskip

\noindent \textbf{(b)} \textit{There is} $\mathfrak{X}$ \textit{such that }$%
\left( \mathcal{M},\mathfrak{X}\right) \models \mathrm{GB}+\mathrm{\Delta }%
_{1}^{1}$\textrm{-}$\mathrm{CA+\Sigma }_{1}^{1}$-$\mathrm{AC}$.\medskip

\noindent \textbf{Proof.~} Suppose $\mathcal{M}$ is a nonstandard model of $%
\mathrm{ZFC}$. By Theorem 3.1 $(b)\Rightarrow (a)$ holds, so we will focus
on establishing $(a)\Rightarrow (b).$ This will be done in two
stages.\medskip

\noindent \textbf{Stage 1.~}We use forcing with set choice functions (as in
Proposition 2.11) to expand $\mathcal{M}$ to a model $(\mathcal{M},f)$ that
satisfies the following properties:\medskip

\noindent (1) $f$ is a global choice function over $\mathcal{M}$, and $(%
\mathcal{M},f)\models \mathrm{ZF}(\overline{\mathbb{L}}_{\mathcal{M}})$%
.\medskip

\noindent (2) $\left( \mathcal{M},f\right) $ is $W$-saturated.\medskip

\noindent Part (b) of Proposition 2.11 assures us that (1) holds. The
verification of (2) involves a careful choice of the generic global choice
function. For this purpose we first verify Lemmas 4.2 and 4.3 below. In
Lemma 4.2 the expression \textquotedblleft $\alpha $ is a Beth-fixed point%
\textit{\textquotedblright } means that $\alpha =\beth (\alpha ),$ where $%
\beth $ is the Beth function. It is well-known that $\mathrm{\alpha }$
is a Beth-fixed point iff $\mathrm{V}(\alpha )$ is a $\Sigma _{1}$%
-elementary submodel of the universe $\mathrm{V}$ of sets.\medskip

\noindent \textbf{4.2.~Lemma.~}(ZFC) \textit{If }$\mathrm{\alpha }$ \textit{%
is a Beth-fixed point and} $\alpha $\textit{\ has countable cofinality, and }%
$\mathcal{N}:=\left( \mathrm{V}(\alpha ),\in \right) $, \textit{then there
is an }$\mathbb{L}_{\mathcal{N}}$-\textit{generic global choice function} $f$
\textit{over} $\mathcal{N}.$ \medskip

\noindent \textbf{Proof.~}\textbf{\ }This is a minor variant of part(a) of
Proposition 2.11 (Forcing Global Choice).\hfill $\square $ (Lemma 4.2)\medskip

\noindent \textbf{4.3.~Lemma.~}\textit{If} $\mathcal{M}$ \textit{is a} $W$-%
\textit{saturated model of }$\mathrm{ZF}(\mathbb{L}_{\mathcal{M}})$, \textit{%
then }$\mathcal{M(\alpha )}\prec _{\mathbb{L}_{\mathcal{M}}}\mathcal{M}$
\textit{for an unbounded collection of }$\alpha \in \mathrm{Ord}^{\mathcal{M}%
}$ \textit{such that} $\mathcal{M}\models \mathrm{cf}(\alpha )=\omega .$
\medskip

\noindent \textbf{Proof.~}\textbf{\ }Fix a nonstandard $\delta \in \mathrm{Ord%
}^{\mathcal{M}}$ and consider the type $p(x,\delta )$ (where $\delta $ is
treated as a parameter) consisting of the formula

\begin{center}
$\left( \delta \in x\right) \wedge \left( x\in \mathrm{Ord}\right) \wedge
\left( \mathrm{cf}(x)=\omega \right) $,
\end{center}

\noindent together with formulae of the form $\mathrm{Ref}_{\varphi }(x)$ as
in Theorem 2.7 (Reflection), where $\varphi $ ranges in $\mathbb{L}_{%
\mathcal{M}}$. It is easy to see that $p(x,y)$ satisfies conditions $(m1)$
and $(m2)$ of part $(m)$ of Definition 2.1. Moreover, by Proposition 2.7 (Reflection), $p(x,\delta )$ also satisfies condition $(m3)$ of the same
definition (since each closed and unbounded subset of ordinals has
unboundedly many members of countable cofinality). Therefore by the
assumption of $W$-saturation of $\mathcal{M}$, $p(x,\mathbb{\delta })$ is
realized in $\mathcal{M}$ by some $\gamma $, which makes it clear that $%
\gamma $ is nonstandard and $\mathcal{M}({\gamma })\prec _{\mathbb{L}_{%
\mathcal{M}}}\mathcal{M}$.\hfill $\square $ (Lemma 4.3)\medskip

\noindent By Proposition 2.7 (Reflection) we can fix a sequence $\left\langle
\alpha _{n}:n<\omega \right\rangle $ that is cofinal in \textrm{Ord}$^{%
\mathcal{M}}$ such that $\mathcal{M(\alpha }_{n}\mathcal{)}\prec _{\mathbb{L}%
_{\mathcal{M}}}\mathcal{M}$ and\textit{\ }$\mathcal{M}\models \mathrm{cf}%
(\alpha _{n})=\omega .$ Then we build an $\mathbb{L}_{\mathcal{M}}$-generic
choice function $f$ over $\mathcal{M}$ by recursively building a sequence of
conditions $\left\langle p_{n}:n<\omega \right\rangle $, as we shall
explain. Thanks to Lemma 4.2 (applied within $\mathcal{M}$) we can get hold
of a condition $p_{1}$ whose domain is $M\mathcal{(\alpha }_{1})$ such that $%
p_{1}$ is $\mathbb{L}_{\mathcal{M(\alpha }_{1})}$-generic over $\mathcal{%
M(\alpha }_{1})$. Generally, given a condition $p_{n}$ in $\mathcal{M}$
whose domain is $M\mathcal{(\alpha }_{n})$ and which is $\mathbb{L}_{%
\mathcal{M(\alpha }_{n})}$-generic over $\mathcal{M(\alpha }_{n})$, we can
use Lemma 4.2 to extend $p_{n}$ to a condition $p_{n+1}$ whose domain is $M%
\mathcal{(\alpha }_{n+1})$, and which is $\mathbb{L}_{\mathcal{M(\alpha }%
_{n+1})}$-generic over $\mathcal{M(\alpha }_{n+1}\mathcal{)}$. Then by the
choice of $\left\langle \alpha _{n}:n<\omega \right\rangle $, the union $f$
of these conditions $\left\langle p_{n}:n<\omega \right\rangle $ will be $%
\mathbb{L}_{\mathcal{M}}$-generic over $\mathcal{M}$. Moreover, $f$ will
have the key property that $f\upharpoonright M(\alpha _{n})$ is $\mathbb{L}_{%
\mathcal{M(\alpha }_{n})}$-generic over $\mathcal{M(\alpha }_{n}\mathcal{)}$
for every $n<\omega $ (and thus truth-and-forcing holds for each of these
approximations). Then thanks again to truth-and-forcing, together with the
fact that $\mathcal{M(\alpha }_{n}\mathcal{)}\prec _{\mathbb{L}_{\mathcal{M}%
}}\mathcal{M}$ for each $n<\omega $, we can conclude:\medskip

\noindent $(\ast )$ $\left( \mathcal{M(\alpha }_{n}\mathcal{)}%
,f\upharpoonright M(\alpha _{n})\right) \prec _{\overline{\mathbb{L}}_{%
\mathcal{M}}}\left( \mathcal{M},f\right) $ for each $n<\omega $. \medskip

\noindent More explicitly, suppose $\left( \mathcal{M(\alpha }_{n}\mathcal{)},f\upharpoonright M(\alpha _{n})\right) \models \varphi(a)$ for some $\overline{\mathbb{L}}_{%
\mathcal{M}}$-formula $\varphi (x)$ and some $a \in M(\alpha _{n})$. Then for some condition $p\in f\upharpoonright M(\alpha _{n})$, we have $\mathcal{M}(\alpha_{n})\models [p \Vdash \varphi(a)]$, and thus by elementarity $\mathcal{M}\models [p \Vdash \varphi(a)]$, which by genericity of $f$ assures us that $\left( \mathcal{M},f\right)\models \varphi(a)$. Note that $(\ast )$ guarantees that the adjunction of the global
choice function $f$ to $\mathcal{M}$ preserves $W$-saturation and concludes
Stage 1 of the proof.\medskip

\noindent \textbf{Stage 2.~}Let $f$ be the global choice function constructed in Stage 1, and let $\mathfrak{X}=\mathrm{Def}_{\overline{\mathbb{L}}_{\mathcal{M}}}(\mathcal{M},f)$, i.e., the family of subsets of $M$ that are definable in $(\mathcal{M},f)$ by some
$\overline{\mathbb{L}}_{\mathcal{M}}$-formula (parameters allowed).

\medskip

\noindent We will treat $f$ as a binary predicate so that variables are the
only terms in $\overline{\mathbb{L}}_{\mathcal{M}}$ (this will slightly
simplify matters in the argument below). By part (b) of Proposition 2.11, $%
\left( \mathcal{M},f\right) \models \mathrm{ZF}(\overline{\mathbb{L}}_{%
\mathcal{M}})$, which makes it clear that $\mathrm{GBC}$ holds in $\left(
\mathcal{M},{\mathfrak{X}}\right) .$ Recall from part (a) of Remark 2.10
that $\mathrm{\Delta }_{1}^{1}$-$\mathrm{CA}$ is provable in $\mathrm{%
GB+\Sigma }_{1}^{1}$-$\mathrm{AC}$, and that in the presence of $\mathrm{GBC}
$, $\mathrm{\Sigma }_{1}^{1}$-$\mathrm{AC}$ is equivalent to $\mathrm{\Sigma
}_{1}^{1}$-$\mathrm{Coll}$. Hence in light of the fact that $\mathrm{GBC}$
holds in $\left( \mathcal{M},{\mathfrak{X}}\right) $ the proof of (b)\ will
be complete once we verify that $\mathrm{\Sigma }_{1}^{1}$-$\mathrm{Coll}$
holds in $\left( \mathcal{M},{\mathfrak{X}}\right) $. For this purpose,
suppose for some parameter $A\in {\mathfrak{X}}$ we have\medskip

\noindent (1) $\left( \mathcal{M},{\mathfrak{X}}\right) \models \forall x\
\exists X\ \psi (x,X,A).$\medskip

\noindent Let $\alpha (m,v)$ be the $\overline{\mathbb{L}}_{\mathcal{M}}$%
-formula that defines $A$, where $m\in M$ is a set parameter. Then\medskip

\noindent (2) $\left( \mathcal{M},{\mathfrak{X}}\right) \models \forall x\
\theta (x),$ where

\begin{center}
$\theta (x):=\bigvee\limits_{\varphi (y,v)\in \overline{\mathbb{L}}_{%
\mathcal{M}}}\exists y$\ $\psi (x,X/\varphi (y,v),A/\alpha (m,v)),$
\end{center}

\noindent and $\psi (X/\varphi (y,v),A/\alpha (m,v),x)$ is the result of
replacing all occurrences of subformulae of the form $w\in X$ (where $w$ is
a variable) in $\psi $ by $\varphi (w,v),$ and replacing all occurrences of
subformulae $w\in A$ in $\psi $ by $\alpha (w,v)$. In these replacements, we
will assume that some variables will be renamed to avoid unintended clashes.
\medskip

\noindent Let $\overline{\mathrm{D}}^{\mathcal{M}}(\alpha )$ consist of all formulae of $%
\overline{\mathbb{L}}_{\infty ,\omega }$ of depth less than $\alpha$ that appear in $\mathcal{M}$. We
claim that (3) below holds.\medskip

\noindent (3) There is some $\alpha \in \mathrm{o}(\mathcal{M})$ such that $%
\mathcal{M}\models \forall x\ \theta _{\alpha }(x)$, where

\begin{center}
$\theta _{\alpha }(x):=\bigvee\limits_{\varphi (y,v)\in \overline{\mathrm{D}}^{\mathcal{M}}
(\alpha )}\exists y$\ $\psi (x,X/\varphi
(y,v),A/\alpha (m,v)).$
\end{center}

\noindent Suppose (3) is false, then we have:\medskip

\noindent (4) $\mathcal{M}\models \exists x\ \lnot \theta _{\alpha }(x)$ for
each $\alpha \in \mathrm{o}(\mathcal{M})$.\medskip

\noindent Consider the $\mathbb{L}_{\mathcal{M}}$-type $p(x):=\left\{ \lnot
\theta _{\alpha }(x):\alpha \in \mathrm{o}(\mathcal{M})\right\} .$ It is
easy to see that $p(x)\in \mathrm{Cod}_{W}(\mathcal{M}).$ By the assumption that (3) is false, for each $%
\alpha \in \mathrm{o}(\mathcal{M})$, $p(x)\cap M(\alpha )$ is realized in $%
\mathcal{M}$, so by $W$-saturation of $\mathcal{M}$, $p(x)$\ is realized in $%
\mathcal{M}$, i.e., $\mathcal{M}\models \exists x\ \lnot \theta (x),$ which
contradicts (2) and completes the verification of (3).

\medskip

Let $B=\mathrm{Def}_{\overline{D}^{\mathcal{M}}(\alpha )}(\mathcal{M},f),$
i.e., the subfamily of $\mathfrak{X}$ consisting of subsets of $M$ that are
definable in $(\mathcal{M},f)$ by some $\overline{\mathbb{L}}_{\mathcal{M}}$%
-formula of depth less than $\alpha$. Note that $B\in \mathfrak{X}$ since there is some $\beta \in
\mathrm{o}(\mathcal{M})$ with $\beta >\mathrm{Depth}(\sigma )$ for each $\sigma \in \overline{\mathrm{D}}^{\mathcal{M}}(\alpha ),$ and \textrm{Sat}$_{\beta }^{%
(\mathcal{M},f)}\in \mathfrak{X}$ by (a minor variant of) Proposition 2.5. Therefore, by (3) we have\medskip

\noindent (5) $\left( \mathcal{M},{\mathfrak{X}}\right) \models \forall x\
\exists y\ \psi (x,X,(B)_{y}).$\medskip

\noindent By quantifying out $B$, (5) readily yields\medskip

\noindent (6) $\left( \mathcal{M},{\mathfrak{X}}\right) \models \exists Y\
\forall x\ \exists y\ \psi (x,(Y)_{y}).$\medskip

\noindent This concludes the verification of $\Sigma _{1}^{1}$-Collection in
$\left( \mathcal{M},{\mathfrak{X}}\right) $.\hfill $\square $ (Theorem
4.1)\medskip

\noindent \textbf{4.4.~Remark.~}The proof of $(b)\Rightarrow (a)$ of Theorem
4.1 does not invoke the countability of $\mathcal{M}$, but the direction $%
(a)\Rightarrow (b)$ does, and indeed this direction of the theorem can fail
for an uncountable model $\mathcal{M}$, e.g., if $\mathcal{M}$ is a
recursively saturated rather classless model of $\mathrm{ZFC+\forall }x%
\mathrm{(V\neq HOD(}x\mathrm{))}$, where $\mathrm{HOD}(x)$ is the class of
sets that are hereditarily ordinal definable from the parameter $x$. More
explicitly, it is well-known that $\mathrm{ZFC+\forall }x\mathrm{(V\neq HOD(}%
x\mathrm{))}$ is consistent, assuming that ZF is consistent.\footnote{%
Easton proved (in his unpublished dissertation \cite{Easton Thesis}) that,
assuming Con($\mathrm{ZF}$), there is a model $\mathcal{M}$ of $\mathrm{ZFC}$
which carries no $\mathcal{M}$-definable global choice function for the
class of pairs in $\mathcal{M}$; and in particular $\exists x\left( \mathrm{V%
}=\mathrm{HOD}(x)\right) $ fails in $\mathcal{M}$. Easton's theorem was
exposited by Felgner \cite[p.231]{Felgner}; for a more recent and
streamlined account, see Hamkins' MathOverflow answer \cite{Joel-failure of
class choice for pairs}.} On the other hand, Kaufmann \cite{Kaufmann}
showed, using the combinatorial principle $\Diamond _{\omega _{1}}$, that
every countable model $\mathcal{M}_{0}$ of ZF has an elementary end
extension $\mathcal{M}$ that is recursively saturated and rather classless,
and later Shelah \cite{Shelah} used an absoluteness argument to eliminate $%
\Diamond _{\omega _{1}}$. Here the rather classlessness of $\mathcal{M}$
means that if $X$ is a subset of $M$ that is piecewise coded in $\mathcal{M}$%
, then $X$ is parametrically definable in $\mathcal{M}$ ($X$ is piecewise
coded in $\mathcal{M}$ means that for every $\alpha \in \mathrm{Ord}^{%
\mathcal{M}},$ $\mathrm{V}^{\mathcal{M}}(\alpha )\cap X$ is coded by an
element of $\mathcal{M}$), then $X$ is parametrically definable in $\mathcal{%
M}$. Therefore if $\mathcal{M}$ is a recursively saturated rather classless
model of $\mathrm{ZFC+\forall }x\mathrm{(V\neq HOD(}x\mathrm{))}$, then by
recursive saturation of $\mathcal{M}$, $\mathcal{M}$ satisfies condition (a)
of Theorem 4.1, but it does not satisfy condition (b) of Theorem 4.1 since
if $\mathcal{M}$ expands to a model $(\mathcal{M},\mathfrak{X})\ $of $%
\mathrm{GB}+\Sigma _{1}^{1}$-AC, then as pointed out in part (b) of Remark
2.10, there is a global choice function $F$ coded in $\mathfrak{X}$. But the
veracity of \textrm{GB} in $(\mathcal{M},\mathfrak{X})$ implies that $F$ is
piecewise coded in $\mathcal{M}$ and therefore $F$ is parametrically
definable in $\mathcal{M}$, which contradicts the fact that $\mathrm{\forall
}x\mathrm{(V\neq HOD(}x\mathrm{))}$ holds in $\mathcal{M}$.

\noindent \textsc{Department of Philosophy, Linguistics, and the Theory of
Science \newline
\noindent University of Gothenburg, Gothenburg, Sweden}\newline
\noindent \texttt{email: ali.enayat@gu.se}

\end{document}